\documentclass{amsart}

\usepackage{amsmath}
\usepackage{amssymb}
\usepackage{commath}
\usepackage{mathtools}
\usepackage[all]{xy}

\usepackage{hyperref}

\theoremstyle{plain}
\newtheorem{theorem}{Theorem}
\newtheorem{lemma}[theorem]{Lemma}
\newtheorem{proposition}[theorem]{Proposition}
\newtheorem{corollary}[theorem]{Corollary}
\theoremstyle{definition}

\newtheorem{remark}[theorem]{Remark}

\numberwithin{theorem}{section}
\numberwithin{equation}{section}

\newcommand{\DI}[1]{\mathrm{D}^\mathrm{I}_{#1 \times #1}}
\newcommand{\UU}[1]{{\mathrm{U}(#1) \times \mathrm{U}(#1)}}
\newcommand{\Pos}[1]{\mathrm{Pos}(n,\mathbb{C})}
\newcommand{\MC}[1]{M_{#1 \times #1}(\C)}

\newcommand{\B}{\mathbb{B}}
\newcommand{\C}{\mathbb{C}}
\newcommand{\N}{\mathbb{N}}
\newcommand{\R}{\mathbb{R}}
\newcommand{\T}{\mathbb{T}}
\newcommand{\Z}{\mathbb{Z}}
\newcommand{\D}{\mathbb{D}}

\newcommand{\GL}{\mathrm{GL}}
\newcommand{\SL}{\mathrm{SL}}
\newcommand{\U}{\mathrm{U}}
\newcommand{\SU}{\mathrm{SU}}
\newcommand{\Spe}{\mathrm{S}}
\newcommand{\Ad}{\mathrm{Ad}}
\newcommand{\ad}{\mathrm{ad}}

\newcommand{\tr}{\mathrm{tr}}

\newcommand{\fg}{\mathfrak{g}}
\newcommand{\fh}{\mathfrak{h}}

\newcommand{\fu}{\mathfrak{u}}
\newcommand{\fsl}{\mathfrak{sl}}

\newcommand{\gl}{\mathfrak{gl}}

\newcommand{\cT}{\mathcal{T}}
\newcommand{\cA}{\mathcal{A}}
\newcommand{\cB}{\mathcal{B}}
\newcommand{\cP}{\mathcal{P}}
\newcommand{\cS}{\mathcal{S}}
\newcommand{\HH}{\mathcal{H}}
\newcommand{\KK}{\mathcal{K}}
\newcommand{\XX}{\mathfrak{X}}

\newcommand{\End}{\mathrm{End}}

\begin{document}

\title[Radial Toeplitz on type I]{Radial-like Toeplitz operators on Cartan domains of type I}

\author{Ra\'ul Quiroga-Barranco}
\address{Centro de Investigaci\'on en Matem\'aticas, Guanajuato, M\'exico}
\email{quiroga@cimat.mx}

\subjclass[2010]{Primary 47B35; Secondary 22D10}

\keywords{Toeplitz operators, representation theory, radial symbols}

\maketitle

\begin{abstract}
	Let $\mathrm{D}^\mathrm{I}_{n \times n}$ be the Cartan domain of type I which consists of the complex $n \times n$ matrices $Z$ that satisfy $Z^*Z < I_n$. For a symbol $a \in L^\infty(\mathrm{D}^\mathrm{I}_{n \times n})$ we consider three radial-like type conditions: 1) left (right) $\U(n)$-invariant symbols, which can be defined by the condition $a(Z) = a\big((Z^*Z)^\frac{1}{2}\big)$ ($a(Z) = a\big((ZZ^*)^\frac{1}{2}\big)$, respectively), and 2) $\mathrm{U}(n) \times \mathrm{U}(n)$-invariant symbols, which are defined by the condition $a(A^{-1}ZB) = a(Z)$ for every $A, B \in \mathrm{U}(n)$. We prove that, for $n \geq 2$, these yield different sets of symbols. If $a$ satisfies 1), either left or right, and $b$ satisfies 2), then we prove that the corresponding Toeplitz operators $T_a$ and $T_b$ commute on every weighted Bergman space.	
	Furthermore, among those satisfying condition 1), either left or right, there exist, for $n \geq 2$, symbols $a$ whose corresponding Toeplitz operators $T_a$ are non-normal. We use these facts to prove the existence, for $n \geq 2$, of commutative Banach non-$C^*$ algebras generated by Toeplitz operators.
\end{abstract}


\section{Introduction}
The theory of Toeplitz operators acting on Bergman spaces has proven to be a very interesting and active line of research. A very general, and at the same time accessible, setup is given by choosing a circular bounded symmetric domain $D$ where we can consider Bergman spaces $\cA^2_\lambda(D)$ and Toeplitz operators acting on them. This can be done for every $\lambda > p-1$, where $p$ is the genus of $D$ (see \cite{UpmeierToepBook}).

An approach to the study of Toeplitz operators is achieved by considering the algebras that they generate when the symbols are restricted to some particular sets. For example, for some subset $\cS \subset L^\infty(D)$, one studies the features of the unital Banach algebra generated by the Toeplitz operators, acting on $\cA^2_\lambda(D)$, whose symbols belong to $\cS$. From now on, we will denote by $\cT^{(\lambda)}(\cS)$ such unital Banach algebra. It turns out that, in many cases, the right choice of set of symbols $\cS$ yields commutative algebras, or at least some setup where studying the commutativity of Toeplitz operators is quite interesting. We refer to \cite{AppuhamyLe2016,AxlerCuckovicRao2000,BauerChoeKoo2015,BVQuasiEllipticI,ChoeKooLee2004,CuckovicLouhichi2008,Le2017,Rodriguez2020,ZorbRadial2002} for just a few examples of this fact.

For a circular bounded symmetric domain $D$, the group of biholomorphism $G$ acts transitively and provides a tool to select special sets of symbols. More precisely, for any closed subgroup $H$ of $G$, we can consider the set of essentially bounded symbols that are $H$-invariant. We will denote by $L^\infty(D)^H$ the space of such $H$-invariant symbols. This has allowed to prove the existence of a variety of commutative $C^*$-algebras of the form $\cT^{(\lambda)}(L^\infty(D)^H)$. Such algebras have shown to be very interesting and complicated in some cases. We refer to \cite{DOQJFA,DQRadial,QVUnitBall1,QVUnitBall2} for some general constructions of commutative $C^*$-algebras obtained with these techniques for the $n$-dimensional unit ball $\B^n$ as well as general irreducible bounded symmetric domains.

Although the general setup of arbitrary bounded symmetric domains is quite important, the case of radial symbols on the unit disk $\D$ can still be considered as particularly interesting. We recall that $a \in L^\infty(\D)$ is called radial if and only if $a(z) = a(|z|)$, for every $z \in \D$. Besides the references mentioned above, we can also refer to \cite{BVquasir-quasih,EsmeralMaximenko2016,GKVRadial,GMVRadial,KorenblumZhu1995,QRadial,QSpseudoconvex,SuarezRadial2008,VasQuasiRadial,VasPseudoH} for examples of research on Toeplitz operators with radial symbols and some generalizations to the unit ball $\B^n$. From these references we would like to highlight \cite{KorenblumZhu1995} where it was discovered the importance of Toeplitz operators with radial symbols. It was first proved there that for any radial symbol $a$ on the unit disk $\D$, the Toeplitz operator $T_a$ acting on the (weightless) Bergman space $\cA^2(\D)$ can be diagonalized with respect to the natural monomial basis. In other words, such Toeplitz operator $T_a$ with radial symbol $a$ preserves the Hilbert direct sum
\begin{equation}\label{eq:intro-Hilbertsum-Korenblum-Zhu}
	\cA^2(\D) = \bigoplus_{n \in \N} \C z^n,
\end{equation}
thus providing its diagonal form. In particular, any Toeplitz operator with radial symbol is normal. Furthermore, it follows that the $C^*$-algebra generated by Toeplitz operators with radial symbols is commutative (see \cite{GKVRadialDisk,GQVJFA}).

We note that a symbol $a \in L^\infty(\D)$ is radial if and only if $a(tz) = a(z)$ for every $t \in \T$ and $z \in \D$. In other words, the space of radial symbols on the unit disk is given by $L^\infty(\D)^\T$, for the $\T$-action on $\D$ by rotations around the origin. This action realizes the isotropy at the origin of the biholomorphisms of $\D$.

A natural problem is to generalize to higher dimensions the notion of radial symbols and to study the features of the corresponding Toeplitz operators. This has already been considered for the $n$-dimensional unit ball $\B^n$. In particular, radial symbols were studied in \cite{GKVRadial}, separately radial symbols (also known as quasi-elliptic symbols) have been studied in \cite{QRadial,QVUnitBall1}, and in other references several variations of these have been considered, for example, \cite{EsmeralMaximenko2016,QSpseudoconvex,VasQuasiRadial,VasPseudoH}. In all these previous works, the corresponding Toeplitz operators have also been studied and used to construct commutative Banach and $C^*$-algebras.

Our main goal here is to consider an alternative generalization, to higher dimensions, of the notion of radial symbol from the unit disk $\D$. More specifically, for the Cartan domains of type I we introduce three possible generalizations of the notion of radial symbols in $\D$. We study the corresponding Toeplitz operators and obtain some commutative Banach algebras that they generate.

We recall that, for every $n \geq 1$, the Cartan domain $\DI{n}$ of type I is the domain of matrices $Z \in M_{n\times n}(\C)$ that satisfy $Z^*Z < I_n$. We note that this defining condition is equivalent to $ZZ^* < I_n$. Also observe that for $n = 1$ we recover the unit disk $\D$.

Given the obvious similarity in the definitions of $\D$ and $\DI{n}$, one might guess that a natural generalization of a radial symbol $a$ for $\DI{n}$ can be given by the condition
\begin{equation}\label{eq:intro-Z*Z}
	a(Z) = a\big((Z^*Z)^\frac{1}{2}\big)
\end{equation}
for every $Z \in \DI{n}$. Alternatively, we can also consider the condition
\begin{equation}\label{eq:intro-ZZ*}
	a(Z) = a\big((ZZ^*)^\frac{1}{2}\big)
\end{equation}
for every $Z \in \DI{n}$. 

On the other hand, since the radial symbols in $\D$ are those invariant by the biholomorphism subgroup fixing the origin, we can consider the corresponding condition for $\DI{n}$. As noted in subsection~\ref{subsec:biholo}, the biholomorphisms fixing the origin in $\DI{n}$ are realized by the linear action of $\UU{n}$ given by
\[
	(A,B)\cdot Z = AZB^{-1}
\]
where $A,B \in \U(n)$ and $Z \in \DI{n}$. Hence, we can consider a third condition for a symbol $a$ in $\DI{n}$, which is
\begin{equation}\label{eq:intro-UUn-invariant}
	a(A^{-1}ZB) = a(Z)
\end{equation}
for every $A,B \in \U(n)$ and $Z \in \DI{n}$. In some sense, we can consider these three conditions as yielding radial-like symbols on the domain $\DI{n}$. Note that, for $n = 1$, these three conditions are clearly equivalent and correspond to the actual radial symbols on $\D$.

However, it turns out that, for $n \geq 2$, the three conditions \eqref{eq:intro-Z*Z}, \eqref{eq:intro-ZZ*} and \eqref{eq:intro-UUn-invariant} are non-equivalent by pairs. To prove this, we introduce in Proposition~\ref{prop:UnUnRn-coordinates} and Corollary~\ref{cor:UnUn(0,1)-coordinates} a sort of polar coordinates similar (but not the same) to those considered in Section~3.4 from \cite{HuaHarmonic} (see Remark~\ref{rmk:polarcoordinates}). Then, Proposition~\ref{prop:UnUn-invariantsymbols} provides a polar coordinates description of the $\UU{n}$-invariant symbols, i.e.~those satisfying \eqref{eq:intro-UUn-invariant}. On the other hand, we prove in Proposition~\ref{prop:leftright-Un-invariant} and Corollary~\ref{cor:leftright-Un-invariant} (see also Remark~\ref{rmk:leftright-Un-invariant}) that conditions \eqref{eq:intro-Z*Z} and \eqref{eq:intro-ZZ*} are equivalent to requiring the symbol $a$ to be invariant under the action of the factors $\U(n)_L = \U(n) \times \{I_n\}$ and $\U(n)_R = \{I_n\} \times \U(n)$, respectively. For this reason, we call the symbols satisfying \eqref{eq:intro-Z*Z} left $\U(n)$-invariant, and we call those satisfying \eqref{eq:intro-ZZ*} right $\U(n)$-invariant (see subsection~\ref{subsec:leftright-Un-symbols}). We also provide a polar coordinates description of both of them in Proposition~\ref{prop:leftright-Un-invariant-varphi}. These results allow us to show that, for $n \geq 2$, all three of the above conditions on a symbol $a$ are different (see Remarks~\ref{rmk:leftright-diagrams} and \ref{rmk:leftright-Un-invariant-varphi}). This will actually bring some richness that allows to construct special Toeplitz operators and Banach algebras.

Hence, our next step is to study the Toeplitz operators whose symbols satisfy one of the three conditions mentioned above. Again, we take as model to generalize the case of the unit disk $\D$. More precisely, we proceed to generalize the Hilbert direct sum decomposition \eqref{eq:intro-Hilbertsum-Korenblum-Zhu}. For this, we note that the latter is precisely the decomposition into irreducible subspaces for a natural $\T$-action on the Bergman space $\cA^2(\D)$.

We consider in subsection~\ref{subsec:UnitaryBergman} a natural unitary representation $\pi_\lambda$ of $\UU{n}$ on the weighted Bergman space $\cA^2_\lambda(\DI{n})$. Then, we obtain in Theorem~\ref{thm:UnUn-isotypic} a Hilbert direct sum which is $\UU{n}$-invariant and that reduces to \eqref{eq:intro-Hilbertsum-Korenblum-Zhu} for $n = 1$ in the weightless case. There is a corresponding result for the actions of $\U(n)_L$ and of $\U(n)_R$ which is obtained in Theorem~\ref{thm:LUnRUn-isotypic}. In all these cases we provide a precise description of the terms in the Hilbert direct sums as modules over the corresponding groups. This turns out to be fundamental since for all three cases the Hilbert direct sum is exactly the same, and can only be distinguished by the unitary actions on the terms of the corresponding groups. In fact, it is the notion of isotypic decomposition from representation theory that allows us to distinguish between the three actions (see Section~\ref{sec:invToeplitz} for further details and definitions).

Furthermore, we prove in Proposition~\ref{prop:a-Ta-Hinvariant} a criterion for a $C^*$-algebra of the form $\cT^{(\lambda)}(L^\infty(\DI{n})^H)$, where $H \subset \UU{n}$ is a closed subgroup, to be commutative. This result is applied by testing whether or not an isotypic decomposition is multiplicity-free (see Section~\ref{sec:invToeplitz} for definitions and further details). This allows us to prove in Theorem~\ref{thm:comm-UnUn} that the $C^*$-algebra $\cT^{(\lambda)}(L^\infty(\DI{n})^{\UU{n}})$ is commutative. Such result was already obtained in \cite{DOQJFA,DQRadial}. However, as noted in Remark~\ref{rmk:UnUn-commutative} our proof provides more information that can be used to find explicit diagonal forms of the corresponding Toeplitz operators.

In contrast with the case of $\UU{n}$-invariance, the left and right $\U(n)$-invariant symbols yield non-commutative $C^*$-algebras. This is proved in Theorem~\ref{thm:noncomm-LUnRUn}, where we show that  $\cT^{(\lambda)}(L^\infty(\DI{n})^{\U(n)_L})$ and $\cT^{(\lambda)}(L^\infty(\DI{n})^{\U(n)_R})$ are not commutative. Nevertheless, we are able to prove in Theorem~\ref{thm:centralizing-LUnRUn} that these algebras centralize each other (see also Corollary~\ref{cor:vonNeumanncentralizing-LUnRUn}). This is a very interesting phenomenon: conditions \eqref{eq:intro-Z*Z} and \eqref{eq:intro-ZZ*}, which generalize radial symbols on $\D$, yield symbols whose Toeplitz operators respectively generate two non-commutative (for $n \geq 2$) $C^*$-algebras that centralize each other. In particular, some degree of commutativity is still present for $n \geq 2$.

These results are possible due to the detailed information provided by the isotypic decompositions obtained for the three groups considered. It also allow us to prove in Theorem~\ref{thm:notnormal-LUnRUn}, for $n \geq 2$, the existence of symbols $a$ satisfying either \eqref{eq:intro-Z*Z} or \eqref{eq:intro-ZZ*} for which the corresponding Toeplitz operators $T^{(\lambda)}_a$ are not normal. This is again in contrast with the behavior observed in the case of the unit disk, where every radial symbol always yields normal Toeplitz operators.

This variety of symbols whose Toeplitz operators are non-normal lead to the existence of some interesting Banach algebras generated by Toeplitz operators. The main results in this regard are Theorems~\ref{thm:noncomm-LUnRUn-oneoperator} and \ref{thm:BanachCommNotC*}, and Corollary~\ref{cor:noncomm-LUnRUn-twooperators} (see also Remark~\ref{rmk:otherBanachcommutative}). In all cases by choosing suitable left and right $\U(n)$-invariant symbols, one obtains, in the case $n \geq 2$, Banach non-$C^*$ algebras generated by Toeplitz operators.

The main tool from representation theory that we used is based on the relationship that exist between the representations of $\U(n)$ and $\GL(n,\C)$. Such relationship comes from the fact that the latter is the complexification of the former (see Lemma~\ref{lem:unitarytrick} and the remarks that follow). From this, a detailed analysis of the well known representation of $\GL(n,\C) \times \GL(n,\C)$ on polynomials over $M_{n \times n}(\C)$, reviewed in Section~\ref{sec:invToeplitz}, provides the means to compute the required isotypic decompositions of Bergman spaces.

The author would like to respectfully dedicate this work to the memory of J\"org Eschmeier, with whom he had the fortune to share an interest on Toeplitz operators acting on higher dimensional domains.

\section{Cartan domains of type I}\label{sec:TypeI}
As observed in the introduction, we will denote by $\DI{n}$ the domain of complex $n \times n$ matrices $Z$ that satisfy $Z^* Z < I_n$. We recall that this condition is equivalent to $Z Z^* < I_n$ (see \cite{HuaHarmonic}). It is well known that $\DI{n}$ is an irreducible circled bounded symmetric domain. The most elementary case is given by $\DI{1} = \D$, the unit disk in the complex plane~$\C$.

The domain $\DI{n}$ is clearly $n^2$-dimensional. Furthermore, this domain has rank $n$, genus $2n$ and its characteristic multiplicities are $a = 2$ and $b = 0$. In particular, it has an unbounded tube-type realization. We refer to \cite{UpmeierToepBook} for the definitions and proofs of these claims.

Let us consider the pseudo-unitary group of matrices $M \in \GL(2n,\C)$ that satisfy the condition $M^\top I_{n,n} \overline{M} = I_{n,n}$ where $I_{n,n}$ is the block diagonal $2n \times 2n$ matrix $I_{n,n} = \mathrm{diag}(I_n, -I_n)$. The group of such matrices $M$ will be denoted by $\U(n,n)$, and it yields the linear isometries of the Hermitian form on $\C^{2n}$ defined by the matrix $I_{n,n}$.

\subsection{The biholomorphism group of $\DI{n}$}
\label{subsec:biholo}
It is well known that the biholomorphism group of $\DI{n}$ is realized by the following action of $\U(n,n)$
\begin{align*}
	\U(n,n) \times \DI{n} &\rightarrow \DI{n} \\
	\begin{pmatrix}
		A & B \\
		C & D
	\end{pmatrix}\cdot Z &= (AZ + B)(CZ + D)^{-1},
\end{align*}
where $A,B,C,D$ have size $n \times n$. We refer to \cite{KnappLieBeyond,MokHermitian} for the details of this fact. Furthermore, this action is transitive and the isotropy subgroup that fixes the origin is the subgroup $\UU{n}$ block diagonally embedded in $\U(n,n)$. In other words, we have the injective homomorphism of Lie groups
\begin{align*}
	\UU{n} &\rightarrow \U(n,n) \\
	(A,B) &\mapsto
			\begin{pmatrix}
				A & 0 \\
				0 & B
			\end{pmatrix},
\end{align*}
that we will consider from now on. In particular, the group of biholomorphisms of $\DI{n}$ that fix the origin is given by the linear action
\begin{align*}
	\UU{n} \times \DI{n} &\rightarrow \DI{n} \\
	(A,B)\cdot Z &= AZB^{-1}.
\end{align*}
Hence, we have a representation of $\DI{n}$ as a homogeneous space given by
\[
	\DI{n} \simeq \U(n,n) / \UU{n}.
\]
However, we note that this representation comes from an action that is not effective. We recall that an action is called effective if the identity element is the only one acting trivially. In this case, it is easy to see that, for every $t \in \T$, the element $(tI_n, tI_n) \in \UU{n}$ acts trivially on $\DI{n}$. One way to avoid this situation is by considering matrices with determinant $1$. More precisely, we consider the groups
\begin{align*}
	\SU(n,n) &= \U(n,n) \cap \SL(2n,\C) \\
	\Spe(\UU{n}) &= \{ (A,B) \mid \det(A)\det(B) = 1 \},
\end{align*}
with the corresponding actions on $\DI{n}$. These new groups continue to realize the whole biholomorphism group and the isotropy subgroup at the origin, respectively, and thus yield the representation 
\[
	\DI{n} \simeq \SU(n,n) / \Spe(\UU{n}).
\]
Even in this case there are elements of $\SU(n,n)$ that act trivially, but now they form a finite group, which is enough to deal with most situations.

We will be mainly interested in the biholomorphisms of $\DI{n}$ that fix the origin, which we just saw that come from the actions of either of the groups $\UU{n}$ or $\Spe(\UU{n})$. For our purposes, it will be more useful to consider the former subgroup rather than the latter. One reason is that this will make more natural to work with two actions of the group $\U(n)$. More precisely, from now on we will refer to the left and the right $\U(n)$-actions which are given respectively by
\begin{equation}\label{eq:leftright-Unactions}
	Z \mapsto AZ, \quad Z \mapsto ZB^{-1},
\end{equation}
for $A, B \in \U(n)$ and $Z \in \DI{n}$. To simplify our notation we will denote
\begin{equation}\label{eq:LUnRUn}
	\U(n)_L = \U(n) \times \{I_n\}, \quad 
		\U(n)_R = \{I_n\} \times \U(n),
\end{equation}
and so the left and right $\U(n)$-actions are given by the actions of $\U(n)_L$ and $\U(n)_R$, respectively. Note the actions of the three groups $\UU{n}$, $\U(n)_L$ and $\U(n)_R$ extend naturally to linear actions on $\MC{n}$, which we will use latter on.

For $n = 1$, we have $\U(1) = \T$, and the left and right $\U(1)$-actions yield the same collection of biholomorphisms on the unit disk $\D$: the rotations around the origin. However, for $n \geq 2$, the group $\U(n)$ is reductive with $1$-dimensional center and its semisimple part is $\SU(n)$ which is in fact simple. The left and right $\U(n)$-actions have in common only the transformations of the form
\[
	Z \mapsto t Z
\]
where $t \in \T$. In other words, if $A, B \in \U(n)$ satisfy $AZ = ZB$, for all $Z \in \DI{n}$, then $A = B = tI_n$ for some $t \in \T$. Since $\DI{n}$ is open in $\MC{n}$, this claim is an easy linear algebra exercise. In particular, for $n \geq 2$, the left and right $\U(n)$-actions on $\DI{n}$ are not the same.

In the rest of this work, the actions of either $\UU{n}$, $\U(n)_L$ or $\U(n)_R$ will be considered both on $\DI{n}$ and $\MC{n}$. As noted above, such actions are linear on $\MC{n}$.

\subsection{Bergman spaces and Toeplitz operators on $\DI{n}$} \label{subsec:BergmanToeplitz}
Let us denote by $\dif v(Z)$ the Lebesgue measure on $\MC{n}$ normalized by the condition $v(\DI{n}) = 1$. The (weightless) Bergman space on $\DI{n}$ is the subspace of holomorphic functions on $\DI{n}$ that belong to $L^2(\MC{n},v)$, and such subspace will be denoted by $\cA^2(\DI{n})$. It is very well known that $\cA^2(\DI{n})$ is a closed subspace and a reproducing kernel Hilbert space whose kernel is the function $K : \DI{n} \times \DI{n} \rightarrow \C$ given by 
\[
	K(Z,W) = \det(I_n - ZW^*)^{-2n}.
\]
We refer to \cite{HuaHarmonic,UpmeierToepBook} for the details of this fact as well as for the claims made in this subsection.

For every $\lambda > 2n-1$, we consider the measure given by
\[
	\det(I_n - ZZ^*)^{\lambda - 2n} \dif v(Z),
\]
which is in fact finite on $\DI{n}$. We let $c_\lambda > 0$ be the constant such that the measure
\[	
	\dif v_\lambda(Z) = c_\lambda 
		\det(I_n - ZZ^*)^{\lambda - 2n} \dif v(Z),
\]
satisfies $v_\lambda(\DI{n}) = 1$. Then, the weighted Bergman space corresponding to a given weight $\lambda > 2n -1$ is the subspace of holomorphic functions that belong to $L^2(\DI{n},v_\lambda)$. We denote this subspace with $\cA^2_\lambda(\DI{n})$ which is, as before, closed and a reproducing kernel Hilbert space. In this case the kernel is the function $K_\lambda : \DI{n} \times \DI{n} \rightarrow \C$ given by
\[
	K_\lambda(Z,W) = \det(I_n - Z W^*)^{-\lambda}.
\]
Note that for $\lambda = 2n$ we recover the weightless Bergman space.

The Bergman projection associated to the Bergman space $\cA^2_\lambda(\DI{n})$ is the projection $B_\lambda : L^2(\DI{n}, v_\lambda) \rightarrow \cA^2_\lambda(\DI{n})$ given by
\[
	B_\lambda(f)(Z)  = 
		\int_{\DI{n}} f(W) K_\lambda(Z, W) \dif v_\lambda(W),
\]
for every $Z \in \DI{n}$.

The above remarks lead to the so-called Toeplitz operators. More precisely, for a function $a \in L^\infty(\DI{n})$, called a symbol, the Toeplitz operator with symbol $a$ is the bounded operator $T_a^{(\lambda)} = T_a \in \cB(\cA^2_\lambda(\DI{n}))$ given by
\begin{align*}
	T_a^{(\lambda)}(f)(Z) 
		= B_\lambda(af)(Z)
		&= \int_{\DI{n}} a(W) f(W) K_\lambda(Z,W) \dif v_\lambda(W) \\
		&= c_\lambda \int_{\DI{n}} 
			\frac{a(W) f(W) \det(I_n - WW^*)^{\lambda-2n} \dif v(W)}{\det(I_n - ZW^*)^{\lambda}},
\end{align*}
for every $Z \in \DI{n}$.

\subsection{A unitary action on Bergman spaces} 
\label{subsec:UnitaryBergman}
The normalized Lebesgue measure $\dif v(Z)$ introduced before is obtained from a corresponding Hermitian inner product on $\MC{n}$. Up to our normalization, such inner product is given by
\[
	(Z,W) \mapsto \tr(ZW^*),
\]
which is clearly invariant under the action of $\UU{n}$ on $\MC{n}$. It follows that the $\UU{n}$-action preserves the normalized Lebesgue measure $\dif v(Z)$. Furthermore, since the function $Z \mapsto \det(I_n - ZZ^*)^{\lambda - 2n}$, defined on $\DI{n}$, is clearly $\UU{n}$-invariant, it follows that the measure $\dif v_\lambda(Z)$ is $\UU{n}$-invariant as well.

From the previous discussion it follows that, for every $\lambda > 2n-1$, we have a unitary representation given by
\begin{align*}
	\UU{n} \times L^2(\DI{n},v_\lambda) &\rightarrow L^2(\DI{n},v_\lambda) \\
	((A,B)\cdot f)(Z) &= f(A^{-1}ZB),
\end{align*}
which is continuous in the strong operator topology. Furthermore, this action clearly leaves invariant the Bergman space $\cA^2_\lambda(\DI{n})$. We will denote by $\pi_\lambda$ this unitary representation of $\UU{n}$ on $\cA^2_\lambda(\DI{n})$.

We note that corresponding properties hold for the subgroups $\U(n)_L$ and $\U(n)_R$. More precisely, the restrictions $\pi_\lambda|_{\U(n)_L}$ and $\pi_\lambda|_{\U(n)_R}$ are unitary representations as well. In words, the left and right $\U(n)$-actions define unitary representations on each weighted Bergman space.

\section{Invariant symbols}\label{sec:invsymbols}
Given the unitary action of $\UU{n}$ on $\cA^2_\lambda(\DI{n})$ introduced in subsection~\ref{subsec:UnitaryBergman} it is useful to consider symbols that are invariant under this group and its subgroups $\U(n)_L$ and $\U(n)_R$. This will yield special types of Toeplitz operators. 

\subsection{Symbols invariant under $\UU{n}$}
\label{subsec:UnxUn-symbols}
We note that there is a natural action of the group $\UU{n}$ on the symbols $a \in L^\infty(\DI{n})$ given by
\[
	((A,B)\cdot a)(Z) = a(A^{-1}ZB)
\]
for every $(A,B) \in \UU{n}$ and for every $Z \in \DI{n}$. A symbol $a$ will be called $\UU{n}$-invariant if it is invariant under this action, in other words if it satisfies
\[
	(A,B)\cdot a = a
\]
for every $(A,B) \in \UU{n}$. We will denote by $L^\infty(\DI{n})^{\UU{n}}$ the space of all $\UU{n}$-invariant symbols. More generally, for every closed subgroup $H \subset \UU{n}$, we will denote by $L^\infty(\DI{n})^H$ the space of symbols $a \in L^\infty(\DI{n})$ that satisfy
\[
h \cdot a = a,
\]
for every $h \in H$. These are also called $H$-invariant symbols.

For every $z \in \C^n$ we will denote by $D(z)$ the diagonal matrix whose diagonal entries are given by the components of $z$. In particular, the collection of matrices $D(t)$ with $t \in \T^n$ is precisely the subgroup of diagonal matrices of $\U(n)$. The assignment $t \mapsto D(t)$ yields an isomorphism between $\T^n$ and the diagonal subgroup of $\U(n)$. From now on we will identify these groups through such assignment. 

From linear algebra (see also \cite{HuaHarmonic}) it follows that for every $Z \in \MC{n}$ there exist $U, V \in \U(n)$ and $x \in \R^n$ satisfying $x_1 \geq \dots \geq x_n \geq 0$ such that
\begin{equation}\label{eq:polarcoordinates}
	Z = U D(x) V.
\end{equation}
This representation is not unique, since for every $t \in \T^n$ we can rewrite it as
\begin{equation}\label{eq:polarcoordinates-Tn}
	Z = U D(x) V = U D(\overline{t}) D(x) D(t) V,
\end{equation}
where we still have $U D(\overline{t}), D(t) V \in \U(n)$. 

In the rest of this work we will denote by $\MC{n}^\times$ the subset of matrices $Z \in \MC{n}$ such that $ZZ^*$ is positive definite with all of its eigenvalues of multiplicity one. We recall that the eigenvalues of $ZZ^*$ and $Z^* Z$ are exactly the same, so we can use $Z^*Z$ instead of $ZZ^*$ in the definition of $\MC{n}^\times$. We also note that $\MC{n}^\times \subset \GL(n,\C)$, and that both of these sets are open conull dense in $\MC{n}$. For $n = 1$, we have $\MC{1}^\times = \GL(1,\C) = \C^\times = \C \setminus \{0\}$, the multiplicative group of non-zero complex numbers, hence our notation.

It is easy to prove that for every $Z \in \MC{n}^\times$ the representation \eqref{eq:polarcoordinates} is unique up to the alternative representations given by \eqref{eq:polarcoordinates-Tn}. This follows, for example, from the computations found in Section~3.4 of \cite{HuaHarmonic}. We can use these remarks to obtain coordinates on $\MC{n}$ almost everywhere that are well-defined up to the ambiguity given by \eqref{eq:polarcoordinates-Tn}.

On the product $\UU{n}$ we consider the $\T^n$-action given~by
\[
	t\cdot (U,V) = (U D(\overline{t}), D(t) V),
\]
where $U, V \in \U(n)$ and $t \in \T^n$. Let us denote by $\overrightarrow{\R}_+^n$ the subset of elements $x$ belonging to $\R^n$ that satisfy $x_1 > \dots > x_n > 0$. We define the smooth map
\begin{align}
	\varphi : \UU{n} \times \overrightarrow{\R}_+^n &\rightarrow \MC{n}^\times
	\label{eq:varphimap} \\
	\varphi(U,V,x) &= U D(x) V. \notag
\end{align}
The next result shows that $\varphi$ parameterizes the open set $\MC{n}^\times$ as the manifold $\UU{n} \times \overrightarrow{\R}_+^n$ up to the $\T^n$-action given above, thus providing coordinates for $\MC{n}$ that are defined almost everywhere and up to such action. Furthermore, it yields an alternative description, up to diffeomorphism, of the set $\MC{n}^\times$ in terms of the groups under consideration. We recall that a submersion is a surjective smooth map whose differential at every point is also surjective. We refer to \cite{KNVol1} for the notion of principal fiber bundle, but we note that the properties established below provide one possible definition.

\begin{proposition}\label{prop:UnUnRn-coordinates}
	The smooth map $\varphi : \UU{n} \times \overrightarrow{\R}_+^n  \rightarrow \MC{n}^\times$ given by \eqref{eq:varphimap} satisfies the following properties.
	\begin{enumerate}
		\item $\varphi$ is a smooth submersion.
		\item The $\T^n$-action on $\UU{n} \times \overrightarrow{\R}_+^n$ is proper and free.
		\item For every $Z \in \MC{n}^\times$, the fiber $\varphi^{-1}(Z)$ is a $\T^n$-orbit.
	\end{enumerate}
	In other words, $\UU{n} \times \overrightarrow{\R}_+^n$ is a principal fiber bundle over $\MC{n}^\times$ with projection $\varphi$ and structure group $\T^n$. In particular, we have an induced diffeomorphism
	\begin{align*}
		\widehat{\varphi} : (\UU{n})/\T^n \times \overrightarrow{\R}_+^n &\rightarrow \MC{n}^\times \\
		\widehat{\varphi}([U,V],x) &= UD(x)V,
	\end{align*}
	where $(\UU{n})/\T^n$ denotes the quotient by the $\T^n$-action on $\UU{n}$ given above, and $[U,V]$ denotes the class of $(U,V) \in \UU{n}$ in such quotient.
\end{proposition}
\begin{proof}
	The properness of the $\T^n$-action is clear from the compactness of this group, and the freeness of the action is a consequence of the freeness of left and right translation actions on groups. This proves (2).
	
	The uniqueness of the representation \eqref{eq:polarcoordinates} up to the alternatives given by \eqref{eq:polarcoordinates-Tn} prove that, for every $Z \in \MC{n}^\times$, the fiber $\varphi^{-1}(Z)$ is precisely a $\T^n$-orbit. This proves (3).

	Clearly the map $\varphi$ is smooth, so it remains to show that its differential at every point is surjective. We note that both the domain and range of $\varphi$ admit natural left and right $\U(n)$-actions which thus define diffeomorphisms. For the domain, these actions are given, respectively, by
	\[
		U\cdot (A,B,x) = (UA,B,x), \quad
			(A,B,x)\cdot V = (A,BV,x).
	\]
	Since $\varphi$ is clearly equivariant for both of these actions, it is enough to prove that $\dif\varphi_{(I_n,I_n,x)}$ is surjective for every $x \in \overrightarrow{\R}_+^n$. A straightforward computation using the multi-linearity of the product of matrices implies that
	\[
		\dif\varphi_{(I_n,I_n,x)}(A,B,v) = 
			AD(x) + D(x)B + D(v),
	\]
	for every $(A,B,v) \in \fu(n) \times \fu(n) \times \R^n$. We will prove that the image of such linear map has real dimension $2n^2 = \dim_\R(\MC{n})$. 
	
	If there exist $A,B \in \fu(n)$ such that $AD(x) + D(x)B = 0$, then applying the adjoint we also have $D(x)A + BD(x) = 0$, because $x \in \R^n$. Hence, for every $j,k$ we have
	\[
		a_{jk} x_k = - x_j b_{jk}, \quad 
			x_j a_{jk} = - b_{jk} x_k,
	\]
	and this clearly yields
	\[
		(x_j^2 - x_k^2) a_{jk} b_{jk} = 0.
	\]
	Since $x_j > x_k > 0$ for every $j > k$, we conclude that both $A,B$ are diagonal and $A = -B$. It follows that the linear map $\fu(n) \times \fu(n) \rightarrow \MC{n}$ given by
	\begin{equation}\label{eq:AD(x)B}
		(A,B) \mapsto AD(x) + D(x)B,
	\end{equation}
	has real $n$-dimensional kernel, and so its image has real dimension $2n^2 -n$.
	
	Now let us assume that there exist $A, B \in \fu(n)$ and $v \in \R^n$ such that
	\[
		AD(x) + D(x)B = D(v).
	\]
	Since both $x,v \in \R^n$, we also have
	\[
		D(x)A + BD(x) = -D(v),
	\]
	and these yield the identities
	\[
		a_{jk} x_k + x_j b_{jk} = \delta_{jk} v_j, \quad
			x_j a_{jk} + b_{jk} x_k = -\delta_{jk} v_j,
	\]
	which for $j\not=k$ give the same equations obtained above. We conclude again that $A = -B$ and that these are diagonal matrices. Hence, for some $y \in \R^n$ we have $A = D(iy)$ and the previous identities reduce to
	\[
		0 = D(iy) D(x) - D(x) D(iy) = D(v),
	\]
	which implies that $v = 0$. Hence, the $n$-dimensional image of the map $v \mapsto D(v)$ is complementary to the image of the map \eqref{eq:AD(x)B}. This proves that the image of the map $\dif \varphi_{(I_n,I_n,x)}$ has real dimension $2n^2$ and so it is surjective.
	
	The rest of the claims now follow from the definition and properties of principal fiber bundles.
\end{proof}

We will use the previous result to introduce coordinates in $\DI{n}$ which, as before, will be defined almost everywhere and up to a $\T^n$-action.

First we state the next result which is an immediate consequence of the condition $Z^*Z < I_n$ that defines $\DI{n}$.

\begin{lemma}\label{lem:DIn(0,1)}
	For a given $Z \in \MC{n}$, and with respect to the representation $Z = U D(x) V$ from \eqref{eq:polarcoordinates} we have: $Z \in \DI{n}$ if and only if $1 > x_1 \geq \dots \geq x_n \geq 0$.
\end{lemma}

We will denote by $\overrightarrow{(0,1)}^n$ the set of all $x \in \R^n$ such that $1 > x_1 > \dots > x_n > 0$. The next result is a consequence of Proposition~\ref{prop:UnUnRn-coordinates} and Lemma~\ref{lem:DIn(0,1)}. We note that the set $\DI{n}\cap \MC{n}^\times$ is an open conull dense subset of the domain $\DI{n}$. For simplicity we will use the same notation $\widehat{\varphi}$ for the map introduced in Proposition~\ref{prop:UnUnRn-coordinates} that now considers the domain $\DI{n}$.

\begin{corollary}\label{cor:UnUn(0,1)-coordinates}
	In the notation of Proposition~\ref{prop:UnUnRn-coordinates} the map given by
	\begin{align*}
		\widehat{\varphi} : (\UU{n})/\T^n \times \overrightarrow{(0,1)}^n &\rightarrow \DI{n}\cap \MC{n}^\times \\
		\widehat{\varphi}([U,V],x) &= U D(x) V,
	\end{align*}
	is a diffeomorphism.
\end{corollary}

\begin{remark}\label{rmk:polarcoordinates}
	It is an easy exercise to verify that the diffeomorphism $\widehat{\varphi}$ considered in Proposition~\ref{prop:UnUnRn-coordinates} and Corollary~\ref{cor:UnUn(0,1)-coordinates} yields, for $n=1$, the usual polar coordinates on $\MC{1}^\times = \C^\times$ and $\DI{1} \cap \MC{1}^\times = \D \setminus \{0\}$, respectively. To see this, we note that there is natural isomorphism
	\[
		(\U(1) \times \U(1))/\T \simeq \T,
	\]
	of Lie groups. Hence, the diffeomorphism $\widehat{\varphi}$ gives, in both cases, a natural generalization of polar coordinates for arbitrary $n \geq 2$. These coordinates are basically equivalent to those considered in Section~3.4 from \cite{HuaHarmonic}. However, both polar coordinates, ours and the ones found in \cite{HuaHarmonic}, use different parameterizing spaces. 
\end{remark}

We now use the coordinates provided by the previous discussion to obtain the next characterization of the $\UU{n}$-invariant symbols.

\begin{proposition}\label{prop:UnUn-invariantsymbols}
	Let $a \in L^\infty(\DI{n})$ be given. Then, the following conditions are equivalent.
	\begin{enumerate}
		\item The symbol $a$ satisfies
			\[
				a(U^{-1}ZV) = a(Z)
			\]
			for every $U,V \in \U(n)$ and $Z \in \DI{n} \cap \MC{n}^\times$.
		\item For every $U, V \in \U(n)$ and $x \in \overrightarrow{(0,1)}^n$ we have $a(UD(x)V) = a(D(x))$.
		\item The function $a\circ \varphi$ defined on $\UU{n} \times \overrightarrow{(0,1)}^n$ depends only on the factor $\overrightarrow{(0,1)}^n$
		\item The function $a\circ \widehat{\varphi}$ defined on $(\UU{n})/\T^n \times \overrightarrow{(0,1)}^n$ depends only on the factor $\overrightarrow{(0,1)}^n$.
	\end{enumerate}
\end{proposition}
\begin{proof}
	From the definitions of $\varphi$ and $\widehat{\varphi}$, we clearly have that (2), (3) and (4) are equivalent. Also, it is immediate that (1) implies~(2).
	
	Let us assume that (2) holds. Let $Z \in \DI{n} \cap \MC{n}^\times$ be given and consider its decomposition $Z = U D(x) V$ as given in \eqref{eq:polarcoordinates}, so that in particular we have $x \in \overrightarrow{(0,1)}^n$. Then, for every $U_1,V_1 \in \U(n)$ we have
	\[
		a(U_1^{-1} Z V_1) = a(U_1^{-1}U D(x) V V_1) 
			= a(D(x)) = a(U D(x) V) = a(Z),
	\]
	and this implies (1).
\end{proof}

\begin{remark}\label{rmk:UnUn-invariant-n=1}
	We recall that a symbol $a \in L^\infty(\D)$ is radial if and only if we have
	\[
		a(z) = a(|z|) = a(|t_1 z t_2|) = a(t_1 z t_2)
	\]
	for all $z \in \D$ and $t_1, t_2 \in \T$, which trivially recovers condition (2) from Proposition~\ref{prop:UnUn-invariantsymbols}. Furthermore, using Remark~\ref{rmk:polarcoordinates} on polar coordinates, the first identity above is also trivially equivalent to (4) of the same proposition. As noted before, in this case $\MC{1}^\times = \C^\times$.
\end{remark}

\subsection{Left and right $\U(n)$-invariant symbols}
\label{subsec:leftright-Un-symbols}
We recall the left and right $\U(n)$-actions on the domain $\DI{n}$ given by \eqref{eq:leftright-Unactions}, which correspond to the actions of the subgroups $\U(n)_L$ and $\U(n)_R$, defined in \eqref{eq:LUnRUn}, respectively. These allow us to consider two additional families of invariant symbols. A symbol $a \in L^\infty(\DI{n})$ will be called left $\U(n)$-invariant if it is invariant under the $\U(n)_L$-action, in other words if we have
\[
	a(UZ) = a(Z),
\]
for every $U \in \U(n)$ and for every $Z \in \DI{n}$. Similarly, the symbol $a \in L^\infty(\DI{n})$ will be called right $\U(n)$-invariant if it is invariant under the $\U(n)_R$-action, which is now equivalent to
\[
	a(ZV) = a(Z),
\]
for every $V \in \U(n)$ and for every $Z \in \DI{n}$. In particular, we can alternatively speak of $\U(n)_L$-invariant symbols and $\U(n)_R$-invariant symbols, respectively. The space consisting of the former will be denote by $L^\infty(\DI{n})^{\U(n)_L}$, while the space consisting of the latter will be denoted by $L^\infty(\DI{n})^{\U(n)_R}$.

Let us denote by $\Pos{n}$ the cone of positive definite Hermitian elements of $\MC{n}$. We recall the following elementary fact from linear algebra.

\begin{lemma}\label{lem:polardecomposition}
	For every $Z \in \GL(n,\C)$ there exist unique elements $U, V \in \U(n)$ and $P, Q \in \Pos{n}$ such that
	\[
		Z = U P = Q V.
	\]
	Furthermore, we have $P = (Z^*Z)^\frac{1}{2}$ and $Q = (ZZ^*)^\frac{1}{2}$.
\end{lemma}

In the notation of Lemma~\ref{lem:polardecomposition} and from now on, we will refer to the expressions $Z = UP$ and $Z = QV$ as the left and right polar decompositions, respectively. The next result provide a characterization of left and right $\U(n)$-invariant symbols in terms of these decompositions. We recall that $\GL(n,\C)$ is open conull dense in $\MC{n}$ and so $\DI{n}  \cap \GL(n,\C)$ has the same properties with respect to $\DI{n}$.

\begin{proposition}\label{prop:leftright-Un-invariant}
	For any symbol $a \in L^\infty(\DI{n})$, the following properties hold.
	\begin{enumerate}
		\item If $a$ satisfies
			\[
				a(Z) = a\big((Z^*Z)^\frac{1}{2}\big)
			\] 
			for every $Z \in \DI{n}$, then $a$ is $\U(n)_L$-invariant. Conversely, if $a$ is $\U(n)_L$-invariant, then we have
			\[
				a(Z) = a\big((Z^*Z)^\frac{1}{2}\big)
			\] 
			for every $Z \in \DI{n}  \cap \GL(n,\C)$.
		\item If $a$ satisfies
			\[
				a(Z) = a\big((ZZ^*)^\frac{1}{2}\big)
			\] 
			for every $Z \in \DI{n}$, then $a$ is $\U(n)_R$-invariant. Conversely, if $a$ is $\U(n)_R$-invariant, then we have
			\[
				a(Z) = a\big((ZZ^*)^\frac{1}{2}\big)
			\] 
			for every $Z \in \DI{n}  \cap \GL(n,\C)$.
	\end{enumerate}
\end{proposition}
\begin{proof}
	If $a$ satisfies the displayed identity in (1) of the statement, then for any given $U \in \U(n)$ we have
	\[
		a(UZ) = a\big(((UZ)^*(UZ))^\frac{1}{2}\big) 
			= a\big((Z^*U^*UZ)^\frac{1}{2}\big)
			= a\big((ZZ^*)^\frac{1}{2}\big) = a(Z)
	\]
	for every $Z \in \DI{n}$, which proves the $\U(n)_L$-invariance of $a$.
	
	If $a$ is $\U(n)_L$-invariant, and we choose $Z \in \DI{n}  \cap \GL(n,\C)$ with left polar decomposition $Z = U (Z^*Z)^\frac{1}{2}$, then we have
	\[
		a(Z) = a\big(U (Z^*Z)^\frac{1}{2}\big) = a\big((Z^*Z)^\frac{1}{2}\big),
	\]
	thus completing the proof of (1). The proof of (2) is similar.
\end{proof}

\begin{remark}\label{rmk:leftright-Un-invariant-n=1}
	For $n =1$, we have $|z| = \sqrt{\overline{z}z} = \sqrt{z\overline{z}}$ for every $z \in \DI{1} = \D$. This yields the trivial fact that, for a symbol $a$, left and right $\T$-invariance are mutually equivalent, and also equivalent to the condition $a(z) = a(|z|)$ for every $z \in \D$. Note than in this case we have $\GL(n,\C) = \C^\times$.
\end{remark}

\begin{remark}\label{rmk:leftright-Un-invariant}
	For a given symbol $a \in L^\infty(\DI{n})$ let us consider the symbol
	\[
		\widehat{a} = \chi_E a.
	\]
	where we take $E = \DI{n} \cap \GL(n,\C)$. In other words, $\widehat{a}$ redefines $a$ as $0$ on the closed null subset of singular matrices in $\DI{n}$. If $Z \in \DI{n}$ is singular, then the matrix $(Z^*Z)^\frac{1}{2}$ is singular as well and so we have $\widehat{a}(Z) = 0 = \widehat{a}\big((Z^*Z)^\frac{1}{2}\big)$. Since the $\U(n)_L$-action leaves invariant the subset $E$, it follows that if $a$ is $\U(n)_L$-invariant, then $\widehat{a}$ is $\U(n)_L$-invariant as well. 
	
	We conclude from these observations and Proposition~\ref{prop:leftright-Un-invariant} that $a$ is $\U(n)_L$-invariant if and only if 
	\[
		\widehat{a}(Z) = \widehat{a}\big((Z^*Z)^\frac{1}{2}\big)
	\]
	for every $Z \in \DI{n}$, and not just for $Z \in \DI{n} \cap \GL(n,\C)$. In other words, any symbol $a \in L^\infty(\DI{n})$ can be redefined outside the $\U(n)_L$-invariant open conull dense subset $\DI{n} \cap \GL(n,\C)$ so that its $\U(n)_L$-invariance is equivalent to the condition
	\[
		a(Z) = a\big((Z^*Z)^\frac{1}{2}\big)
	\]
	to hold for every $Z \in \DI{n}$. On the other hand, we note that for continuous symbols this equivalence is immediate from the density of $\DI{n} \cap \GL(n,\C)$ in $\DI{n}$, in which case there is no need to redefine any of the values of $a$.
	
	We also note that similar remarks hold for $\U(n)_R$-invariant symbols.
\end{remark}

As a consequence of Remark~\ref{rmk:leftright-Un-invariant} we obtain the next result. 

\begin{corollary}\label{cor:leftright-Un-invariant}
	Up to the identification of symbols that are equal almost everywhere, the left and right $\U(n)$-invariant symbols are given as follows.
	\begin{align*}
		L^\infty(\DI{n})^{\U(n)_L} &= 
			\{ a \in L^\infty(\DI{n}) \mid a(Z) = a\big( (Z^*Z)^\frac{1}{2}\big) \text{ for every } Z \in \DI{n} 
			\}, \\
		L^\infty(\DI{n})^{\U(n)_R} &= 
			\{ a \in L^\infty(\DI{n}) \mid a(Z) = a\big( (ZZ^*)^\frac{1}{2}\big) \text{ for every } Z \in \DI{n} \}.		
	\end{align*}
\end{corollary}

Straight from the definition of invariance we have
\[
	L^\infty(\DI{n})^{\UU{n}} 
		= L^\infty(\DI{n})^{\U(n)_L} \cap L^\infty(\DI{n})^{\U(n)_R}.
\]
In particular, the space $L^\infty(\DI{n})^{\UU{n}}$ is a subspace of both $L^\infty(\DI{n})^{\U(n)_L}$ and $L^\infty(\DI{n})^{\U(n)_R}$. As noted in Remark~\ref{rmk:leftright-Un-invariant-n=1}, these three spaces are all the same for $n=1$. However, we will show that for $n \geq 2$ there are plenty of both left and right $\U(n)$-invariant symbols that are not $\UU{n}$-invariant. We achieve this through the next result which corresponds to Proposition~\ref{prop:UnUn-invariantsymbols}. As noted above, we will consider that two given functions are the same if they agree almost everywhere.

\begin{proposition}\label{prop:leftright-Un-invariant-varphi}
	Let $a \in L^\infty(\DI{n})$ be a given symbol. Then, the following conditions are equivalent.
	\begin{enumerate}
		\item The symbol $a$ is left (right, respectively) $\U(n)$-invariant.
		\item The function $([U,V],x) \mapsto a \circ \widehat{\varphi}([U,V],x) = a(UD(x)V)$ defined on the set $(\UU{n})/\T^n \times \overrightarrow{(0,1)}^n$ is independent of $U$ (independent of $V$, respectively).
		\item There is a measurable function $f$ defined on $(\T^n \backslash \U(n)) \times \overrightarrow{(0,1)}^n$ (defined on $(\U(n)/\T^n) \times \overrightarrow{(0,1)}^n$, respectively) such that 
		\[
			a \circ \widehat{\varphi} = f \circ \rho
		\]
		where $\rho$ is given by $([U,V],x) \mapsto ([V],x)$ (given by $([U,V],x) \mapsto ([U],x)$, respectively).
	\end{enumerate}
	In particular, the assignment 
	\[
		f \mapsto f \circ \rho \circ \widehat{\varphi}^{-1}
	\]
	establishes a one-to-one correspondence between either $L^\infty((\T^n \backslash \U(n)) \times \overrightarrow{(0,1)}^n)$ or $L^\infty((\U(n)/T^n) \times \overrightarrow{(0,1)}^n)$, and the $\U(n)_L$-invariant or $\U(n)_R$-invariant symbols, respectively, where $\rho$ is either of the assignments $([U,V],x) \mapsto ([V],x)$ or $([U,V],x) \mapsto ([U],x)$, respectively.
\end{proposition}
\begin{proof}
	We will only consider the case of left $\U(n)$-invariance since the other case is proved similarly.
	
	If $a$ is $\U(n)_L$-invariant, then we have 
	\[
		a \circ \widehat{\varphi}([U,V],x) = a(UD(x)V) 
			= a(D(x)V) = a \circ \widehat{\varphi}([I_n,V],x),
	\]
	and so it follows immediately that (1) implies (2).
	
	Let us assume that (2) holds. For every $V \in \U(n)$ and its corresponding class $[V] \in \U(n)$ we define
	\[
		f([V],x) = a \circ \widehat{\varphi}([I_n,V],x) = a(D(x) V),
	\]
	for every $x \in \overrightarrow{(0,1)}^n$. From the definition of the corresponding quotients it is clear that $[I_n,V] = [I_n,V_1]$ implies $[V] = [V_1]$. Furthermore, if $[V] = [V_1]$, then there exist $t \in \T^n$ such that $V = D(t)V_1$ and we have
	\[
		a(D(x) V) = a(D(x)D(t)V_1) = a(D(t)D(x)V_1) = a(D(x)V_1),
	\]
	where we have used (2) in the last identity. It follows that $f$ is a well defined function on $\T^n \backslash \U(n) \times \overrightarrow{(0,1)}^n$. Also, using any measurable section of $\rho$ it is easy to see that $f$ is measurable. This proves (3).
	
	Let us now assume that (3) holds. Hence, we have
	\[
		a(UD(x)V) = a\circ\widehat{\varphi}([U,V],x) = f([V],x) 
			= a\circ\widehat{\varphi}([I_n,V],x) = a(D(x)V)
	\]
	for every $U,V \in \U(n)$ and $x \in \overrightarrow{(0,1)}^n$. Hence, for every $Z \in \DI{n} \cap \MC{n}^\times$ with decomposition $Z = UD(x)V$ as in \eqref{eq:polarcoordinates} and $U_1 \in \U(n)$ we have
	\[
		a(U_1 Z) = a(U_1U D(x) V) = a(D(x) V) = a(U D(x) V) = a(Z).
	\]
	This proves the $\U(n)_L$-invariance of $a$ on a conull subset of $\DI{n}$ which, according to our convention, is enough to conclude (1).
	
	For the last claim it is enough to recall from Corollary~\ref{cor:UnUn(0,1)-coordinates} that $\widehat{\varphi}$ is a diffeomorphism from $(\UU{n})/\T^n \times \overrightarrow{(0,1)}^n$ onto the open conull dense subset $\DI{n} \cap \MC{n}^\times$ of $\DI{n}$, and so it has a smooth inverse.
\end{proof}

\begin{remark}\label{rmk:leftright-diagrams}
	We note that condition (3) from Proposition~\ref{prop:leftright-Un-invariant-varphi} for the left $\U(n)$-invariant case yields the commutative diagram
	\[
	\xymatrix{
		\DI{n} \ar[r]^a & \C \\
		(\UU{n})/\T^n \times \overrightarrow{(0,1)}^n 
		\ar[r] \ar[u]_{\widehat{\varphi}} &
		(\T^n\backslash\U(n)) \times \overrightarrow{(0,1)}^n 
		\ar[u]^f
	}
	\]
	where the lower horizontal arrow is given by the assignment 
	\[
		([U,V],x) \mapsto ([V],x).
	\]
	And for the right $\U(n)$-invariant case it yields the commutative diagram
	\[
	\xymatrix{
		\DI{n} \ar[r]^a & \C \\
		(\UU{n})/\T^n \times \overrightarrow{(0,1)}^n 
		\ar[r] \ar[u]_{\widehat{\varphi}} &
		(\U(n)/\T^n) \times \overrightarrow{(0,1)}^n 
		\ar[u]^f
	}
	\]
	where the lower horizontal arrow is given by the assignment 
	\[
		([U,V],x) \mapsto ([U],x).
	\]
\end{remark}

\begin{remark}\label{rmk:leftright-Un-invariant-varphi}
	By comparing Propositions~\ref{prop:UnUn-invariantsymbols} and \ref{prop:leftright-Un-invariant-varphi}, we observe that the $\UU{n}$-invariant symbols correspond to functions defined on $\overrightarrow{(0,1)}^n$, while the left and right $\U(n)$-invariant symbols correspond to functions defined on either $(\T^n\backslash\U(n)) \times \overrightarrow{(0,1)}^n$ or $(\U(n)/\T^n) \times \overrightarrow{(0,1)}^n$, respectively. The first factor is in both cases a manifold with real dimension $n(n-1)$. In fact, both manifolds are easily seen to be diffeomorphic through the assignment $[U] \mapsto [U^{-1}]$, so one can work with either once this identification has been considered.
	
	For $n = 1$, this manifold is just a point, which corresponds to the equivalence of the invariance with respect to the three groups involved. However, for $n \geq 2$, the manifold has positive dimension and Proposition~\ref{prop:leftright-Un-invariant-varphi} proves the existence of plenty of either left or right $\U(n)$-invariant symbols that are not $\UU{n}$-invariant. More precisely, there are as many of such symbols as elements in $L^\infty(\T^n\backslash\U(n))$. In particular, there are plenty of symbols that are either left or right $\U(n)$-invariant but not both. 
\end{remark}

\section{Toeplitz operators with invariant symbols}
\label{sec:invToeplitz}
The principal tool to obtain our results will be representation theory. Hence, we recall some notions and refer to \cite{BtD,GW,KnappLieBeyond} for further details and proofs.

If $H$ is a compact Lie group, then a unitary representation of $H$ on a Hilbert space $\HH$ is a strong topology continuous linear action $\pi: H \times \HH \rightarrow \HH$. In this case, we also say that $\HH$ is an $H$-module. For $\KK$ a closed subspace of $\HH$ we will say that $\KK$ is $H$-invariant, or an $H$-submodule, if $\pi(g)(\KK) = \KK$ for every $g \in H$. If $\KK$ is $H$-invariant and its only $H$-invariant subspaces are $0$ and $\KK$ itself, then $\KK$ is called irreducible. It is known that, since $H$ is compact, all the irreducible $H$-modules are finite dimensional and every $H$-module, finite dimensional or not, is a Hilbert direct sum of irreducible $H$-modules. 

If $\HH_1$ and $\HH_2$ are two given $H$-modules, with corresponding representations $\pi_1$ and $\pi_2$, then a bounded operator $T : \HH_1 \rightarrow \HH_2$ such that $T\circ \pi_1(h) = \pi_2(h) \circ T$, for every $h \in H$, is called an intertwining operator. If there is a unitary intertwining operator $\HH_1 \rightarrow \HH_2$, then we say that the $H$-modules are isomorphic. As usual this yields an equivalence relation. For $\HH = \HH_1 = \HH_2$, the algebra of all intertwining operators is denoted by $\End_H(\HH)$, and it is a von Neumann algebra. 

We will denote by $\widehat{H}$ the family of all equivalence classes of irreducible $H$-modules. We now state our main abstract tool from representation theory of compact groups.

\begin{proposition}\label{prop:isotypicabstract}
	Let $H$ be a compact group with a unitary representation on a Hilbert space $\HH$. For every equivalence class $[\KK] \in \widehat{H}$ let us denote by
	\[
		\HH_{[\KK]} = \sum \{ W \subset \HH \mid
			W \text{ is an $H$-submodule in the class of $[\KK]$} \}.
	\]
	Then, the following properties hold.
	\begin{enumerate}
		\item For every $[\KK] \in \widehat{H}$, the subspace $\HH_{[\KK]}$ is a (closed) $H$-submodule of $\HH$. This is called the isotypic component of $\HH$ associated to the class $[\KK]$.
		\item There is a Hilbert direct sum 
		\[
			\HH = \bigoplus_{[\KK] \in \widehat{H}} \HH_{[\KK]}.
		\]
		This is called the isotypic decomposition of $\HH$ as an $H$-module. If $T$ belongs to $\End_H(\HH)$, then $T$ preserves this Hilbert direct sum.
		\item The von Neumann algebra $\End_H(\HH)$ is commutative if and only if for every $[\KK] \in \widehat{H}$ the isotypic component $\HH_{[\KK]}$ is either $0$ or irreducible. If this is the case, then we say that the isotypic decomposition is multiplicity-free, and every $T$ belonging to $\End_H(\HH)$ acts by a constant multiple of the identity on each isotypic component.
	\end{enumerate}
\end{proposition}

\subsection{Invariance with respect to closed subgroups}
\label{subsec:closedinvariance}
We are interested on the unitary representations $\pi_\lambda$ (for $\lambda > 2n-1$) of the group $\UU{n}$ on the weighted Bergman spaces $\cA^2_\lambda(\DI{n})$, as well as on the restrictions $\pi_\lambda|_{\U(n)_L}$ and $\pi_\lambda|_{\U(n)_R}$ (see subsection~\ref{subsec:UnitaryBergman}). However, for simplicity, in this subsection we will consider any closed subgroup $H \subset \UU{n}$. Then there is a correspondence between $H$-invariant symbols and intertwining Toeplitz operators. We recall from Section~\ref{sec:invsymbols} that a symbol $a \in L^\infty(\DI{n})$ is $H$-invariant if and only if $a \circ h = a$ for every $h \in H$, and that the family of all such symbols is denoted by $L^\infty(\DI{n})^H$. We also denote by $\pi_\lambda|_H$ the restriction of the representation $\pi_\lambda$ to $H$. 

We recall that for a set $\cS \subset L^\infty(\DI{n})$ of symbols we denote by $\cT^{(\lambda)}(\cS)$ the unital Banach algebra generated by Toeplitz operators with symbols in $\cS$ and acting on $\cA^2_\lambda(\DI{n})$, where the weight $\lambda > 2n-1$. In particular, for any closed subgroup $H \subset \UU{n}$, the algebra $\cT^{(\lambda)}(L^\infty(\DI{n})^H)$ is a $C^*$-algebra. This is a well known consequence of the fact that $L^\infty(\DI{n})^H$ is conjugation invariant.

The next proposition is the main tool from representation theory that we will use in the setup of Toeplitz operators. Similar results can be found in \cite{DOQJFA}, but we will use the one below that requires a simplified proof which is enough for our purposes. Furthermore, our particular statement will provide some useful tools.

\begin{proposition}\label{prop:a-Ta-Hinvariant}
	If $H \subset \UU{n}$ is a closed subgroup, then the following holds for every $\lambda > 2n-1$.
	\begin{enumerate}
		\item A symbol $a \in L^\infty(\DI{n})$ is $H$-invariant if and only if $T^{(\lambda)}_a$ is intertwining for $\pi_\lambda|_H$. In other words, we have 
		\[
			\cT^{(\lambda)}(L^\infty(\DI{n})^H) \subset \End_H(\cA^2_\lambda(\DI{n})).
		\]
		\item For every $T \in \End_H(\cA^2_\lambda(\DI{n}))$ and any finite dimensional $\pi_\lambda|_H$-invariant subspace $W \subset \cA^2_\lambda(\DI{n})$, there exists a symbol $a \in L^\infty(\DI{n})^H$ such that 
		\[
			\langle T f, g \rangle_\lambda 
				= \langle T^{(\lambda)}_a f, g \rangle_\lambda
		\]
		for every $f,g \in W$. 
		\item If all the isotypic components associated to $\pi_\lambda|_H$ are finite dimensional, then the $C^*$-algebra $\cT^{(\lambda)}(L^\infty(\DI{n})^H)$ is commutative if and only if the isotypic decomposition of the restriction $\pi_\lambda|_H$ is multiplicity-free.
	\end{enumerate}
\end{proposition}
\begin{proof}
	The proof of (1) is an easy exercise that uses the formula
	\[
		\langle T^{(\lambda)}_a f, g \rangle_\lambda = 
			\langle a f, g \rangle_\lambda,
	\]
	for every $f, g \in \cA^2_\lambda(\DI{n})$, together with the fact that $\pi_\lambda$ is unitary.
	
	To prove (2), let $T$ and $W$ be given as in its statement. Theorem~2 from~\cite{EnglisDensity} proves the existence of a symbol $a \in L^\infty(\DI{n})$ such that 
	\[
		\langle T f, g \rangle_\lambda 
		= \langle T^{(\lambda)}_a f, g \rangle_\lambda
	\]
	for every $f,g \in W$. Let us consider the symbol given by
	\[
		\widehat{a}(Z) = \int_H a(h^{-1}\cdot Z) \dif h,
	\]
	where $\dif h$ is the probability Haar measure of $H$. Clearly $\widehat{a}$ belongs to $L^\infty(\DI{n})^H$, as a consequence of the left and right invariance of the Haar measure for compact groups. Then, we can compute as follows for every $f,g \in W$
	\begin{align*}
		\langle T^{(\lambda)}_{\widehat{a}} f, g \rangle_\lambda 
			&= \langle \widehat{a} f, g \rangle_\lambda \\
			&= \int_{\DI{n}} \int_H a(h^{-1}\cdot Z) f(Z) 
					\overline{g(Z)} \dif h \dif v_\lambda(Z) \\
			&= \int_H \int_{\DI{n}} a(Z) f(h \cdot Z) 
				\overline{g(h \cdot Z)} \dif v_\lambda(Z) \dif h \\	
			&= \int_H \langle a \pi_\lambda(h)^{-1} f ,
					\pi_\lambda(h)^{-1} g \rangle_\lambda \dif h \\
			&= \int_H \langle T^{(\lambda)}_a \pi_\lambda(h)^{-1} f ,
				\pi_\lambda(h)^{-1} g \rangle_\lambda \dif h \\
			&= \int_H \langle T \pi_\lambda(h)^{-1} f ,
				\pi_\lambda(h)^{-1} g \rangle_\lambda \dif h \\
			&= \int_H \langle T f , g \rangle_\lambda \dif h 
					= \langle T f , g \rangle_\lambda,
	\end{align*}
	where we have used the properties of $T$, $W$ and $a$. This proves (2) for the symbol~$\widehat{a}$.
	
	For (3) we first observe that (1) implies that the commutativity of the $C^*$-algebra $\cT^{(\lambda)}(L^\infty(\DI{n})^H)$ follows from that of $\End_H(\cA^2_\lambda(\DI{n}))$, which occurs when $\pi_\lambda$ is multiplicity-free. 
	
	Conversely, let us assume that $\pi_\lambda$ is not multiplicity-free. Let us choose a non-irreducible isotypic component $W$ of $\pi_\lambda$ which, by assumption, is finite dimensional. Then, there exist $T_1, T_2 \in \End_H(\cA^2_\lambda(\DI{n}))$ that necessarily preserve $W$ and whose restrictions to $W$ do not commute with each other. By (2), there exist symbols $a_1, a_2 \in L^\infty(\DI{n})^H$ such that
	\[
		\langle T_j f, g \rangle_\lambda 
			= \langle T^{(\lambda)}_{a_j} f, g \rangle_\lambda
	\]
	for every $f,g \in W$, and $j=1,2$. By (1) and Proposition~\ref{prop:isotypicabstract}(2), the Toeplitz operators $T^{(\lambda)}_{a_1}$ and $T^{(\lambda)}_{a_1}$ preserve $W$ as well, because the latter is an isotypic component. It follows that
	\[
		T_1|_W = T^{(\lambda)}_{a_1}|_W, \quad 
			T_2|_W = T^{(\lambda)}_{a_2}|_W,
	\]
	which implies that $T^{(\lambda)}_{a_1}$ and $T^{(\lambda)}_{a_1}$ do not commute with each other and so the $C^*$-algebra $\cT^{(\lambda)}(L^\infty(\DI{n})^H)$ is not commutative.
\end{proof}

Propositions~\ref{prop:isotypicabstract} and \ref{prop:a-Ta-Hinvariant} lead us to find the isotypic decompositions of the $\UU{n}$-action and the left and right $\U(n)$-actions to determine the commutativity of the algebras generated by Toeplitz operators whose symbols are invariant under such actions.

\subsection{Toeplitz operators with $\UU{n}$-invariant symbols} 
\label{subsec:UUn-invariance}
Let us denote by $\cP(\MC{n})$ the vector space of polynomials on $\MC{n}$, and by $\cP^d(\MC{n})$ the subspace of those that are homogeneous of degree $d$. Then, we have an algebraic direct sum
\[
	\cP(\MC{n}) = 
		\bigoplus_{d=0}^\infty \cP^d(\MC{n})
\]
which is invariant with respect to transformations of the form $Z \mapsto tZ$ for every $t \in \T$. Since these transformations belong to the $\UU{n}$-action, they induce a unitary representations on the weighted Bergman spaces. Furthermore, it is well known that $\cP(\MC{n})$ is dense in every weighted Bergman space $\cA^2_\lambda(\DI{n})$ (see~\cite{UpmeierToepBook}). Hence, for every $\lambda > 2n-1$, we have a Hilbert direct sum decomposition
\begin{equation}\label{eq:BergmanPoly}
	\cA^2_\lambda(\DI{n}) = 
		\bigoplus_{d=0}^\infty \cP^d(\MC{n})
\end{equation}
that corresponds to the isotypic decomposition of the $\T$-action just described. More precisely, every subspace $\cP^d(\MC{n})$ is the isotypic component associated to the character $\chi_{-d}(t) = t^{-d}$ of $\T$.

Since the $\UU{n}$-action is linear, it follows that it preserves the Hilbert direct sum \eqref{eq:BergmanPoly}. Hence, to obtain a decomposition of $\cA^2_\lambda(\DI{n})$ into irreducible $\UU{n}$-submodules it is enough to do so for each term of \eqref{eq:BergmanPoly}. We will achieve this by using representation theory for the group $\GL(n,\C)$. 

The next result establishes an equivalence between representations of suitable complex and compact Lie groups, which justify our passage from $\U(n)$ to $\GL(n,\C)$. It is well known from representation theory as part of the so-called Weyl's unitary trick. Hence, we provide a sketch of the proof and refer to \cite{GW, KnappLieBeyond} for further details. 

\begin{lemma}\label{lem:unitarytrick}
	Let $G$ be a connected complex Lie group and $H$ a closed subgroup. Assume that the Lie algebra of $G$ satisfies $\fg = \fh \oplus i\fh$, where $\fh$ is the Lie algebra of $H$. If $\pi : G \rightarrow \GL(W)$ is a finite dimensional complex representation, then the following properties hold.
	\begin{enumerate}
		\item If $W_0 \subset W$ is a (complex) subspace, then $W_0$ is $G$-invariant if and only if it is $H$-invariant.
		\item The representation $\pi$ is irreducible if and only if $\pi|_H$ is irreducible.
	\end{enumerate}
	Furthermore, if $\pi_j : G \rightarrow \GL(W_j)$, for $j=1, 2$, are two given finite dimensional representations, then $W_1 \simeq W_2$ as $G$-modules if and only if $W_1 \simeq W_2$ as $H$-modules.
\end{lemma}
\begin{proof}
	Clearly (1) implies (2), so we will prove the former.
	
	It is obvious that $G$-invariance of a subspace implies its $H$-invariance. So we consider an $H$-invariant subspace $W_0$ and show that it is $G$-invariant.
	
	Let $\dif\pi : \fg \rightarrow \gl(W)$ be the induced representation on Lie algebras obtained by differentiation. In particular, we have a commutative diagram
	\[
		\xymatrix{
			\fg \ar[r]^{\dif\pi} \ar[d] & \gl(W) \ar[d] \\
			G \ar[r]^{\pi} & \GL(W)	
		}
	\]
	where the vertical arrows are the corresponding exponential maps. Since $\pi(h)(W_0) = W_0$ for every $h \in H$, it follows that
	\[
		\pi(\exp(rX))(w) \in W_0
	\]
	for every $X \in \fh$, $w \in W_0$ and $r \in \R$. Hence, differentiation with respect to $r$ yields $\dif\pi(X)(W_0) \subset W_0$ for every $X \in \fh$. Since $W_0$ is complex, $\dif\pi$ is complex linear and $\fg = \fh \oplus i\fh$ this implies that $\dif\pi(X)(W_0) \subset W_0$ for every $X \in \fg$. We conclude that
	\[
		\exp(\dif\pi(X))(w) 
			= \sum_{j=0}^\infty \frac{\dif\pi(X)^j}{j!}(w) \in W_0
	\]
	for every $X \in \fg$ and $w \in W_0$. The above commutative diagram now implies that
	\[
		\pi(\exp(X))(W_0) = \exp(\dif\pi(X))(W_0) = W_0
	\]
	for every $X \in \fg$. Since $G$ is connected it is generated by $\exp(\fg)$ and this yields the $G$-invariance of $W_0$.
	
	To prove the last claim, note that the non-trivial part is showing that $W_1 \simeq W_2$ as $H$-modules implies $W_1 \simeq W_2$ as $G$-modules. Let $T : W_1 \rightarrow W_2$ be an isomorphism of $H$-modules. In particular, we have
	\[
		T\circ \pi_1(h) = \pi_2(h) \circ T,
	\]
	for every $h \in H$. As before, taking $h = \exp(rX)$ and differentiating with respect to $r$ we obtain
	\[
		T \circ \dif \pi_1(X) = \dif\pi_2(X) \circ T,
	\]
	for every $X \in \fh$. But this implies that the same identity holds for every $X \in \fg$ since $T$, $\dif\pi_1$ and $\dif\pi_2$ are complex linear. Using the series expansion of the exponential of linear maps we conclude that
	\[
		T \circ \exp(\dif\pi_1(X)) = \exp(\dif\pi_2(X)) \circ T,
	\]
	for every $X \in \fg$. Commutative diagrams similar to the one used above imply that 
	\[
		T \circ \pi_1(\exp(X)) = \pi_2(\exp(X)) \circ T,
	\]
	for every $X \in \fg$. Since $G$ is connected it is generated by $\exp(\fg)$ and we conclude that $T$ is an isomorphism of $G$-modules.
\end{proof}

Lemma~\ref{lem:unitarytrick} clearly applies to $G = \GL(n,\C)$ and $H = \U(n)$. This will allow us to describe isotypic decompositions involving $\U(n)$ in terms of the irreducible representations of $\GL(n,\C)$. For this reason, we will recall some of the properties of the irreducible rational representations of $\GL(n,\C)$. This will be enough for our purposes. We refer to \cite{GW} for further details, definitions and proofs of the results and claims found in the rest of this subsection.

Let us denote by $\C^{\times n}$ the subgroup of diagonal matrices in $\GL(n,\C)$. In particular, we have $\T^n = \U(n) \cap \C^{\times n}$, the subgroup of diagonal matrices in $\U(n)$. Then, the Lie algebra of $\C^{\times n}$ is $\C^n$ and the Lie algebra of $\T^n$ is $i\R^n$, both considered as spaces of diagonal matrices. In other words, we will use from now on the embeddings
\[
	\C^n \hookrightarrow \gl(n,\C), \quad \C^{\times n} \hookrightarrow \GL(n,\C), 
\]
of Lie algebras and Lie groups, respectively, both given by the assignment $z \mapsto D(z)$. Note that a corresponding remark applies for the Lie group $\U(n)$.

Let $\pi : \GL(n,\C) \rightarrow \GL(W)$ be a rational representation of $\GL(n,\C)$, where $W$ is a finite dimensional complex vector space. Then, for every $\mu \in \C^{n*}$, the complex dual space of $\C^n$, we will denote 
\[
	W(\mu) = \{ w \in W \mid \dif\pi(X)(w) = \mu(X) w,
		\text{ for every } X \in \C^n \}.
\]
If $W(\mu) \not= 0$, then this subspace is called a weight space and the functional $\mu$ is called a weight, both associated to $W$ as $\GL(n,\C)$-module. If we denote by $\XX(W)$ the set of weights for the $\GL(n,\C)$-module $W$, then it is well known (see Section~3.1.3 from \cite{GW}) that $W$ is the direct sum of its weight spaces
\[
	W = \bigoplus_{\mu \in \XX(W)} W(\mu).
\]
A particular, case is given by the adjoint representation of $\GL(n,\C)$
\begin{align*}
	\Ad : \GL(n,\C) &\rightarrow \GL(\gl(n,\C)) \\
	\Ad(g)(X) &= g X g^{-1},
\end{align*}
whose differential is well known to be the adjoint representation of $\gl(n,\C)$
\begin{align*}
	\ad : \gl(n,\C) &\rightarrow \gl(\gl(n,\C)) \\
	\ad(X)(Y) &= [X,Y].	
\end{align*}
The collection of non-zero weights for the adjoint representation of $\GL(n,\C)$ is called the set of roots and it will be denoted by $\Phi$. 

Let $\{e_j\}_{j=1}^n$ be the canonical basis of the dual space $\C^{n*}$ of $\C^n$. In other words, we define $e_j(z) = z_j$, for every $z \in \C^n$. Recall that $\C^n$ is being identified with the subspace of diagonal matrices in $\MC{n}$. Then, it is easily seen that the set of roots for $\GL(n,\C)$ is given by
\[
	\Phi = \{ e_j - e_k \mid j,k = 1, \dots, n,\; j\not= k\}.
\]
We note that the functionals $\{e_j\}_{j=1}^n$ are real valued on $\R^n$. Hence, we can identify its real span with the dual space $\R^{n*}$ of $\R^n$. In particular, we will consider from now on $\Phi \subset \R^{n*}$, and also that $e_j \in \R^{n*}$ for every $j =1, \dots, n$. The lexicographic order with respect to the ordered basis $e_1, \dots, e_n$ yields a partial order on $\R^{n*}$. We will denote by $\Phi_+$ the set of positive roots with respect to this order. In particular, we now have
\[
	\Phi_+ = \{ e_j - e_k \mid 1 \leq j < k \leq n\}.
\]

As before, let us consider any given finite dimensional rational representation $\pi : \GL(n,\C) \rightarrow \GL(W)$. Then, we define the root order in $\XX(W)$ by 
\[
	\mu_1 \prec \mu_2 \Longleftrightarrow \mu_2 - \mu_1 = \alpha_1 + \dots + \alpha_m \text{ for some } \alpha_1, \dots, \alpha_m \in \Phi_+.
\]

We note that the same constructions and properties considered so far apply to the Lie group $\SL(n,\C)$ and its Lie algebra $\fsl(n,\C)$ without any essential change. The only modification that has to be done is to replace the diagonal subgroup $\C^{\times n}$ of $\GL(n,\C)$ by the subgroup $\C^{\times n}_0$ of $z \in \C^{\times n}$ such that $z_1 \cdot \dots \cdot z_n = 1$. Correspondingly, the Lie algebra of $\C^\times_0$ is the subspace $\C^n_0$ of $z \in \C^n$ such that $z_1 + \dots + z_n = 0$. One advantage of considering the subgroup $\SL(n,\C)$ is that it is semisimple (for $n \geq 2$) and we have at our disposal the Theorem of the Highest Weight for Lie algebras (see Section~3.2.1 from \cite{GW}). This allows us to obtain the next result. We sketch the additional arguments required to obtain the needed statement for the Lie group $\SL(n,\C)$ from the results found in \cite{GW}. Note that for a finite dimensional rational representation $\pi$ of either $\GL(n,\C)$ or $\SL(n,\C)$ (the cases under consideration) the theory of weights is obtained from the representation $\dif \pi$ of either $\gl(n,\C)$ or $\fsl(n,\C)$, respectively.

\begin{proposition}\label{prop:SL-highestweight}
	The finite dimensional rational representations of $\SL(n,\C)$ satisfy the following properties.
	\begin{enumerate}
		\item If $\pi : \SL(n,\C) \rightarrow \GL(W)$ is an irreducible finite dimensional rational representation, then there exists a unique weight $\mu_W \in \XX(W)$ such that $\mu \prec \mu_W$ for every $\mu \in \XX(W) \setminus \{\mu_W\}$. The weight $\mu_W$ is called the highest weight of~$W$.
		\item If $\pi_j : \SL(n,\C) \rightarrow \GL(W_j)$, for $j =1,2$, are two irreducible finite dimensional rational representations, then $W_1 \simeq W_2$ as $\SL(n,\C)$-modules if and only if $\mu_{W_1} = \mu_{W_2}$. In other words, two such irreducible representations are equivalent if and only if they have the same highest weight.
	\end{enumerate}
\end{proposition}
\begin{proof}
	The existence of highest weights for the noted irreducible representations is a consequence of Corollary~3.2.3 from \cite{GW}, thus proving (1).
	
	Let us consider two representations $\pi_1$ and $\pi_2$ as in (2). If $W_1$ and $W_2$ are isomorphic as $\SL(n,\C)$-modules, then differentiation proves that they are isomorphic as $\fsl(n,\C)$-modules. From the definitions it is easy to see that this implies $\XX(W_1) = \XX(W_2)$, and so we conclude that $\mu_{W_1} = \mu_{W_2}$.
	
	We now assume that $\mu_{W_1} = \mu_{W_2}$. Then, Theorem~3.2.5 from \cite{GW} implies that $W_1 \simeq W_2$ as $\fsl(n,\C)$-modules. Let $T : W_1 \rightarrow W_2$ be an isomorphism of $\fsl(n,\C)$-modules. Then, we have
	\[
		T \circ \dif\pi_1(X) = \dif\pi_2(X)\circ T,
	\]
	for every $X \in \fsl(n,\C)$. At this point we can repeat the arguments found at the end of the proof of Lemma~\ref{lem:unitarytrick} to show that $T$ is an isomorphism of $\SL(n,\C)$-modules.
\end{proof}

The fundamental dominant weights (see \cite{GW}) for the group $\SL(n,\C)$ are the linear functionals given by the restrictions
\[
	\omega_j = (e_1 + \dots + e_j)|_{\C^n_0}
\]
for every $j =1, \dots, n-1$, which belong to $\C^{n*}_0$, the dual space of $\C^n_0 \subset \fsl(n,\C)$. Note that the restriction of $e_1 + \dots + e_n$ to $\C^n_0$ is $0$.

The set of dominant weights for the group $\SL(n,\C)$ is given by (see Section~3.1.4 from~\cite{GW})
\[
	P_{++}(\SL(n,\C)) = \bigoplus_{j=1}^{n-1} \N \omega_j.
\]
With this notation, the next result is a consequence of Section~3.1.4 and Theorem~5.5.21 from \cite{GW}.

\begin{corollary}\label{cor:SL-highestweight}
	The set $P_{++}(\SL(n,\C))$ is precisely the collection of the highest weights of the irreducible rational representations of the group $\SL(n,\C)$. In particular, there is a one-to-one correspondence between the family of equivalence classes of irreducible rational representations of $\SL(n,\C)$ and the set $P_{++}(\SL(n,\C))$. Such correspondence assigns to every equivalence class the highest weight of any of its elements.
\end{corollary}

\begin{remark}\label{rmk:SL-highestweight}
	For every $\mu \in P_{++}(\SL(n,\C))$ we will denote by $W^\mu$ a $\SL(n,\C)$-module with highest weight $\mu$. In particular, $W^\mu$ is well defined up to an isomorphism of $\SL(n,\C)$-modules. Thus, the inverse of the correspondence stated in Corollary~\ref{cor:SL-highestweight} is given by $\mu \mapsto [W^\mu]$.
\end{remark}

We will now use the previous constructions and follow Section~5.5.4 from \cite{GW} to describe the irreducible rational representations of $\GL(n,\C)$. The set of dominant weights of the group $\GL(n,\C)$ is the set $P_{++}(\GL(n,\C))$ of elements $\mu \in \C^{n*}$, the dual space of $\C^n \subset \gl(n,\C)$, that can be written as
\begin{equation}\label{eq:mu-ej}
	\mu = m_1 e_1 + \dots + m_n e_n, 
\end{equation}
where $m_1 \geq \dots \geq m_n$ and  $m_j \in \Z$ for every $j = 1, \dots, n$.

For every $\mu \in P_{++}(\GL(n,\C))$ given by the expression \eqref{eq:mu-ej} we will consider the element of $P_{++}(\SL(n,\C))$ given by
\begin{equation}\label{eq:mu0-from-mu}
	\mu_0 = (m_1 - m_2)\omega_1 + \dots + (m_{n-1} - m_n) \omega_{n-1}.
\end{equation}

With the previous notation, the next result describes the irreducible rational representations of $\GL(n,\C)$. This is basically a restatement of Theorem~5.5.22 from \cite{GW}.

\begin{proposition}\label{prop:GL-highestweighttheorem}
	There is a one-to-one correspondence between the set of dominant weights $P_{++}(\GL(n,\C))$ and the equivalence classes of irreducible rational representations of $\GL(n,\C)$. This correspondence assigns to every $\mu \in P_{++}(\GL(n,\C))$ an irreducible rational representation $\pi^\mu : \GL(n,\C) \rightarrow \GL(W^\mu)$ that satisfies the following properties, where $\mu$ has the representation given by \eqref{eq:mu-ej}.
	\begin{enumerate}
		\item The restriction $\pi^\mu|_{\SL(n,\C)}$ is an irreducible representation of $\SL(n,\C)$ with highest weight $\mu_0$ given by \eqref{eq:mu0-from-mu}.
		\item The restriction $\pi^\mu|_{\C^\times I_n}$ yields the representation of the diagonal subgroup $\C^\times \simeq \C^\times I_n \subset \GL(n,\C)$ given by the action on $W^\mu$ through the character $z \mapsto z^{m_1 + \dots + m_n}I_n$.
	\end{enumerate}
\end{proposition}

We will use the previous constructions to study suitable representations of $\U(n)$. Hence, the following result will be very useful.

\begin{lemma}\label{lem:Un-highestweighttheorem}
	With the notation of Proposition~\ref{prop:GL-highestweighttheorem}, for every $\mu \in P_{++}(\GL(n,\C))$ the irreducible representation $\pi^\mu : \GL(n,\C) \rightarrow \GL(W^\mu)$ restricted to $\U(n)$ is irreducible as well. In other words, $W^\mu$ is an irreducible $\U(n)$-module for every $\mu \in P_{++}(\GL(n,\C))$. Furthermore, if $\mu_1, \mu_2 \in P_{++}(\GL(n,\C))$ then $W^{\mu_1} \simeq W^{\mu_2}$ as $\U(n)$-modules if and only if $\mu_1 = \mu_2$.
\end{lemma}
\begin{proof}
	Lemma~\ref{lem:unitarytrick}(2) implies the first part of the statement. The second part of the statement follows from Proposition~\ref{prop:GL-highestweighttheorem} and the last claim of Lemma~\ref{lem:unitarytrick}. In both cases, to apply Lemma~\ref{lem:unitarytrick} we take $G = \GL(n,\C)$ and $H = \U(n)$.
\end{proof}

The representation $\pi_\lambda$ of $\UU{n}$ on $\cP(\MC{n})$ has a natural extension to the representation
\begin{align*}
	\GL(n,\C) \times \GL(n,\C) \times \cP(\MC{n}) &\rightarrow 
			\cP(\MC{n}) \\
			(A,B)\cdot p(Z) &= p(A^{-1}ZB),
\end{align*}
which clearly preserves $\cP^d(\MC{n})$ for every $d \in \N$. Furthermore, the induced representation of $\GL(n,\C) \times \GL(n,\C)$ on each of these subspaces is rational (see~\cite{GW}). We will now describe its decomposition into irreducible submodules. First, we recall some properties of representations of a product of groups and tensor products that we will use freely. We refer to \cite{BtD,GW} for further details and proofs.

Recall that for two given Lie groups $H_1, H_2$ and corresponding finite dimensional modules $W_1, W_2$, the tensor product $W_1 \otimes W_2$ admits a natural representation of $H_1 \times H_2$. If both $W_1$ and $W_2$ are irreducible, then $W_1 \otimes W_2$ is irreducible as well. Furthermore, if $W_1, W_1'$ and $W_2,W_2'$ are finite dimensional irreducible modules over $H_1$ and $H_2$, respectively, then $W_1 \otimes W_2$ and $W_1' \otimes W_2'$ are isomorphic over $H_1 \times H_2$ if and only if $W_1 \simeq W_1'$ and $W_2 \simeq W_2'$ over $H_1$ and $H_2$, respectively.

We now recall the definition of a special type of dominant weight that appears in the representation of $\GL(n,\C) \times \GL(n,\C)$ on the space of polynomials given above. A dominant weight $\mu$ will be called non-negative if $m_1 \geq \dots \geq m_n \geq 0$, where $\mu$ is given by \eqref{eq:mu-ej}. In this case we will write
\[
	|\mu| = m_1 + \dots + m_n,
\]
and we call $|\mu|$ the size of $\mu$. Let us denote by $P(n)$ the set of all non-negative dominant weights of $\GL(n,\C)$, and by $P_d(n)$ the subset of those with size $d \in \N$. In particular, we have $P_d(n) \subset P(n) \subset P_{++}(\GL(n,\C))$ for every $d \in \N$. Then, the next result is obtained using Theorem~5.6.7 from~\cite{GW}.

\begin{proposition}\label{prop:GL-poly-isotypic}
	For every $\mu \in P(n)$ there is a $\GL(n,\C) \times \GL(n,\C)$-submodule $\cP^\mu(\MC{n})$ of $\cP(\MC{n})$ such that the following properties are satisfied.
	\begin{enumerate}
		\item For every $d \in \N$ we have a direct sum decomposition
			\[
				\cP^d(\MC{n}) = \bigoplus_{\mu \in P_d(n)} 
						\cP^\mu(\MC{n})
			\]
			which is $\GL(n,\C) \times \GL(n,\C)$-invariant.
		\item For every $\mu \in P(n)$, the spaces $\cP^\mu(\MC{n})$ and	$W^{\mu *} \otimes W^\mu$ are isomorphic as $\GL(n,\C) \times \GL(n,\C)$-modules, where $W^\mu$ is the $\GL(n,\C)$-module given by Proposition~\ref{prop:GL-highestweighttheorem} and $W^{\mu *}$ is its dual $\GL(n,\C)$-module. In particular, the direct sum in (1) is the isotypic decomposition for the $\GL(n,\C) \times \GL(n,\C)$-action and it is multiplicity-free. 
		\item The algebraic direct sum 
			\[
				\cP(\MC{n}) = \bigoplus_{\mu \in P(n)}
							\cP^\mu(\MC{n})
			\]
			is the isotypic decomposition for the $\GL(n,\C) \times \GL(n,\C)$-action, and it is multiplicity-free.
	\end{enumerate}
\end{proposition}
\begin{proof}
	Claims (1) and (2) follow directly from Theorem~5.6.7 from \cite{GW}. The direct sum in (3) and its invariance is clear. On the other hand, for any two given summands $\cP^{\mu_1}(\MC{n})$ and $\cP^{\mu_2}(\MC{n})$ the action of the subgroup
	\[
		\{ (zI_n, z^{-1}I_n) \mid z \in \C^\times \} \simeq \C^\times
	\]
	of $\GL(n,\C) \times \GL(n,\C)$ on them is given by the characters 
	\[
		\chi_{-2|\mu_1|}(z) = z^{-2|\mu_1|}, \quad 
			\chi_{-2|\mu_2|}(z) = z^{-2|\mu_2|},
	\]
	since they are subspaces of $\cP^{|\mu_1|}(\MC{n})$ and $\cP^{|\mu_2|}(\MC{n})$, respectively. It follows from this together with (1) and (2) that the direct sum in (3) consists of mutually non-isomorphic irreducible $\GL(n,\C) \times \GL(n,\C)$-submodules. This completes the proof of (3).	
\end{proof}

As a consequence, we obtain the next result for Bergman spaces.

\begin{theorem}\label{thm:UnUn-isotypic}
	For every $\lambda > 2n -1$, we have a Hilbert direct sum 
	\[
		\cA^2_\lambda(\DI{n}) = \bigoplus_{\mu \in P(n)}
					\cP^\mu(\MC{n}).
	\]
 	which is the isotypic decomposition for the $\UU{n}$-action. Furthermore, for every $\mu \in P(n)$ the subspace $\cP^\mu(\MC{n})$ is isomorphic to $W^{\mu *} \otimes W^\mu$ as a module over $\UU{n}$. In particular, this $\UU{n}$-action is multiplicity-free.
\end{theorem}
\begin{proof}
	Lemma~\ref{lem:unitarytrick}(2) and Proposition~\ref{prop:GL-poly-isotypic}(2) imply that for every $\mu \in P(n)$ the subspace $\cP^\mu(\MC{n})$ is irreducible and isomorphic to $W^{\mu *} \otimes W^\mu$, both over $\UU{n}$. The last claim of Lemma~\ref{lem:unitarytrick} and Proposition~\ref{prop:GL-poly-isotypic}(3) imply that all such subspaces are mutually non-isomorphic over $\UU{n}$. Note that we applied Lemma~\ref{lem:unitarytrick} to the case $G = \GL(n,\C) \times \GL(n,\C)$ and $H = \UU{n}$.
	
	It follows that, for every $d \in \N$, the direct sum 
	\[
	\cP^d(\MC{n}) = \bigoplus_{\mu \in P_d(n)} \cP^\mu(\MC{n})
	\]
	from Proposition~\ref{prop:GL-poly-isotypic}(1), is the isotypic decomposition for the $\UU{n}$-action. In particular, this sum is orthogonal with respect to any $\UU{n}$-invariant inner product.
	
	The previous arguments and \eqref{eq:BergmanPoly} imply that the sum in the statement is indeed a Hilbert direct sum. The rest of the claims of the statement follow as well from the previous remarks.
\end{proof}

Theorem~\ref{thm:UnUn-isotypic} yields the following result. It provides an alternative proof to some of the results from \cite{DQRadial}.

\begin{theorem}\label{thm:comm-UnUn}
	For any $\lambda > 2n -1$, the $C^*$-algebra $\cT^{(\lambda)}(L^\infty(\DI{n})^{\UU{n}})$ acting on $\cA^2_\lambda(\DI{n})$ is commutative. Furthermore, every operator belonging to $\cT^{(\lambda)}(L^\infty(\DI{n})^{\UU{n}})$ preserves the Hilbert direct sum from Theorem~\ref{thm:UnUn-isotypic} and acts by a constant multiple of the identity on each of its summands.
\end{theorem}
\begin{proof}
	By Theorem~\ref{thm:UnUn-isotypic}, the restriction $\pi_\lambda|_{\UU{n}}$ is multiplicity-free and so Proposition~\ref{prop:isotypicabstract}(3) implies that $\End_{\UU{n}}(\cA^2_\lambda(\DI{n}))$ is commutative. Hence, the first claim follows from Proposition~\ref{prop:a-Ta-Hinvariant}(1) applied to $H = \UU{n}$. The second claim follows from Proposition~\ref{prop:isotypicabstract}(3) and Proposition~\ref{prop:a-Ta-Hinvariant}(1).
\end{proof}

\begin{remark}\label{rmk:UnUn-commutative}
	The commutativity claim from Theorem~\ref{thm:comm-UnUn} was proved in \cite{DQRadial} for the domain $\DI{n}$ among others. However, the corresponding isotypic decomposition was not explicitly computed in \cite{DQRadial}. Hence, Theorems~\ref{thm:UnUn-isotypic} and \ref{thm:comm-UnUn} together add information that can be used to study with more detail the operators belonging to the $C^*$-algebra $\cT^{(\lambda)}(L^\infty(\DI{n})^{\UU{n}})$.
\end{remark}

\subsection{Toeplitz operators with left and right $\U(n)$-invariant symbols} 
\label{subsec:leftrightUn-invariance}
Let us now consider the left and right $\U(n)$-actions on the weighted Bergman spaces $\cA^2_\lambda(\DI{n})$ for every $\lambda > 2n-1$. These actions correspond to the subgroups $\U(n)_L$ and $\U(n)_R$, respectively, of $\UU{n}$. In particular, the Hilbert direct sum from Theorem~\ref{thm:UnUn-isotypic} is invariant under both the left and right $\U(n)$-actions. However, as we will show, the terms of such direct sum are no longer irreducible for them.

The next result yields the previous claims by describing the isotypic decompositions for the restrictions $\pi_\lambda|_{\U(n)_L}$ and $\pi_\lambda|_{\U(n)_R}$.

\begin{theorem}\label{thm:LUnRUn-isotypic}
	For every $\lambda > 2n-1$, the Hilbert direct sum
	\[
		\cA^2_\lambda(\DI{n}) = \bigoplus_{\mu \in P(n)}
				\cP^\mu(\MC{n}),
	\]
	is the isotypic decomposition for both the left and right $\U(n)$-actions. More precisely, for every $\mu \in P(n)$, the subspace $\cP^\mu(\MC{n})$ is an isotypic component for both the left and right $\U(n)$-actions. Furthermore, these decompositions satisfy the following properties.
	\begin{enumerate}
		\item For the representation $\pi_\lambda|_{\U(n)_L}$ on $\cA^2_\lambda(\DI{n})$ we have an isomorphism of $\U(n)$-modules
		\[
			\cA^2_\lambda(\DI{n}) = \bigoplus_{\mu \in P(n)}
					\cP^\mu(\MC{n})
				\simeq \bigoplus_{\mu \in P(n)} 
					\bigoplus_{j=1}^{\dim W^\mu} W^{\mu *},
		\]
		obtained from the fact that, for every $\mu \in P(n)$, the subspace $\cP^\mu(\MC{n})$ is isomorphic, as $\U(n)_L$-module, to the sum of $\dim W^\mu$ copies of the $\U(n)$-module~$W^{\mu *}$.
		\item For the representation $\pi_\lambda|_{\U(n)_R}$ on $\cA^2_\lambda(\DI{n})$ we have an isomorphism of $\U(n)$-modules
		\[
			\cA^2_\lambda(\DI{n}) = \bigoplus_{\mu \in P(n)}
					\cP^\mu(\MC{n})
				\simeq \bigoplus_{\mu \in P(n)} 
					\bigoplus_{j=1}^{\dim W^\mu} W^\mu,
		\]
		obtained from the fact that, for every $\mu \in P(n)$, the subspace $\cP^\mu(\MC{n})$ is isomorphic, as $\U(n)_R$-module, to the sum of $\dim W^\mu$ copies of the $\U(n)$-module~$W^\mu$.
	\end{enumerate}
	In particular, for $n \geq 2$, the restricted representations $\pi_\lambda|_{\U(n)_L}$ and $\pi_\lambda|_{\U(n)_R}$ are not multiplicity-free.
\end{theorem}
\begin{proof}
	The Hilbert direct sum holds by Theorem~\ref{thm:UnUn-isotypic} and it is left and right $\U(n)$-invariant as noted above. To prove the rest of the statement we will only consider the left $\U(n)$-action since the case of the right $\U(n)$-action can be handled similarly.
	
	Let $\mu \in P(n)$ be given. By Theorem~\ref{thm:UnUn-isotypic} there is an isomorphism 
	\[
		T : \cP^\mu(\MC{n}) \rightarrow W^{\mu *} \otimes W^\mu
	\]	
	of $\UU{n}$-modules. In particular, $T$ is an isomorphism of $\U(n)_L$-modules. Note that the $\U(n)_L$-module structure on the target is given by
	\[
		U \cdot (w^* \otimes w) = (U \cdot w^*) \otimes w,
	\]
	for every $U \in \U(n)$, $w^* \in W^{\mu *}$ and $w \in W^\mu$. Hence, for every $w \in W^\mu \setminus \{0\}$ the subspace $W^{\mu *} \otimes w$ is a $\U(n)_L$-submodule naturally isomorphic to the $\U(n)$-module $W^{\mu *}$. If we choose a basis $w_1, \dots, w_{\dim W^\mu}$ of $W^\mu$, then we conclude that
	\[
		W^{\mu *} \otimes W^\mu 
			= \bigoplus_{j=1}^{\dim W^\mu} W^{\mu *} \otimes w_j
			\simeq \bigoplus_{j=1}^{\dim W^\mu} W^{\mu *}
	\]
	is a decomposition into irreducible $\U(n)_L$-submodules, where the indicated isomorphism is over $\U(n)$. These remarks and $T$ yield the isomorphism of $\U(n)$-modules
	\[
		\cP^\mu(\MC{n}) \simeq \bigoplus_{j=1}^{\dim W^\mu} W^{\mu *}.
	\]
	On the other hand, for every $\mu_1, \mu_2 \in P(n)$, we have $W^{\mu_1 *} \simeq W^{\mu_2 *}$ if and only $W^{\mu_1} \simeq W^{\mu_2}$, where both isomorphisms are considered over $\U(n)$. By Lemma~\ref{lem:unitarytrick} (for $G = \GL(n,\C)$ and $H = \U(n)$) and Proposition~\ref{prop:GL-highestweighttheorem} such conditions are also equivalent to $\mu_1 = \mu_2$.
	
	The previous arguments show that the subspace $\cP^\mu(\MC{n})$, where $\mu \in P(n)$, is an isotypic component for the $\U(n)_L$-action. They also prove the claims in~(1).
	
	Finally, let us assume that $n \geq 2$. Then, it is well known that the group $\U(n)$ has infinitely many irreducible representations of the form $W^\mu$, for some $\mu \in P(n)$, and so that $\dim W^\mu \geq 2$. For example, if $d \in \N$, then we have $d e_1 \in P(n)$ and $W^{d e_1} \simeq \cP^d(\C^n)^*$ as $\U(n)$-modules, and so $\dim W^{d e_1} \geq 2$ for $d \geq 1$. Since $\cP^\mu(\MC{n})$ is the sum of $\dim W^\mu$ irreducible $\U(n)$-modules, we conclude that the isotypic decomposition for $\pi_\lambda|_{\U(n)_L}$ is not multiplicity-free.
\end{proof}

We now obtain the next result for Toeplitz operators whose symbols are left or right $\U(n)$-invariant but not necessarily both.

\begin{theorem}\label{thm:noncomm-LUnRUn}
	For every $n \geq 2$ and $\lambda > 2n-1$, both of the unital $C^*$-algebras $\cT^{(\lambda)}(L^\infty(\DI{n})^{\U(n)_L})$ and $\cT^{(\lambda)}(L^\infty(\DI{n})^{\U(n)_R})$ are not commutative. However, every operator belonging to either of these algebras preserve the Hilbert direct sum from Theorem~\ref{thm:LUnRUn-isotypic}.
\end{theorem}
\begin{proof}
	We only consider $\cT^{(\lambda)}(L^\infty(\DI{n})^{\U(n)_L})$, since the other algebra can be treated similarly. Note that we are assuming that $n \geq 2$ and $\lambda > 2n -1$.
	
	Theorem~\ref{thm:LUnRUn-isotypic} shows that the restriction $\pi_\lambda|_{\U(n)_L}$ is not multiplicity-free and that its isotypic components are finite dimensional. Thus, it follows from Proposition~\ref{prop:a-Ta-Hinvariant}(3) that the $C^*$-algebra $\cT^{(\lambda)}(L^\infty(\DI{n})^{\U(n)_L})$ is not commutative. The second claim follows from Proposition~\ref{prop:a-Ta-Hinvariant}(1) and Proposition~\ref{prop:isotypicabstract}(2) applied to the subgroup $H = \U(n)_L$.
\end{proof}

We now show that the non-commutativity of the $C^*$-algebras from Theorem~\ref{thm:noncomm-LUnRUn} is quite strong. It will also be useful in constructing commutative Banach algebras generated by Toeplitz operators that are non-$C^*$.

\begin{theorem}\label{thm:notnormal-LUnRUn}
	Let us assume that $n \geq 2$. Then, for every $\lambda > 2n-1$, there exist $a \in L^\infty(\DI{n})^{\U(n)_L}$ and $b \in L^\infty(\DI{n})^{\U(n)_R}$ such that the Toeplitz operators $T^{(\lambda)}_a$ and $T^{(\lambda)}_b$ are not normal.
\end{theorem}
\begin{proof}
	We will consider only the case of the subgroup $\U(n)_L$, since the case of $\U(n)_R$ is proved similarly.
	
	By Theorem~\ref{thm:LUnRUn-isotypic}, the isotypic decomposition for $\pi_\lambda|_{\U(n)_L}$ is not multiplicity-free. Hence, there exist $\mu \in P(n)$ such that the subspace $\cP^\mu(\MC{n})$ is not irreducible over $\U(n)_L$. Furthermore, by the claims in Theorem~\ref{thm:LUnRUn-isotypic}(1) it follows that
	\[
		\End_{\U(n)_L}(\cP^\mu(\MC{n})) \simeq \MC{m}
	\]
	as $C^*$-algebras, where $m = \dim W^\mu \geq 2$. Let $A \in \MC{m}$ be some matrix which is not normal. This yields an operator $T_A \in \End_{\U(n)_L}(\cP^\mu(\MC{n}))$ which is not normal either. Extend $T_A$ by zero in the rest of the terms of the Hilbert direct sum from Theorem~\ref{thm:LUnRUn-isotypic} to obtain an operator $T \in \End_{\U(n)_L}(\cA^2(\DI{n}))$ such that
	\[
		T|_{\cP^\mu(\DI{n})} = T_A.
	\]
	By Proposition~\ref{prop:a-Ta-Hinvariant}(2) there exist $a \in L^\infty(\DI{n})^{\U(n)_L}$ such that
	\[
		\langle T f, g \rangle_\lambda 
			= \langle T^{(\lambda)}_a f, g \rangle_\lambda
	\]
	for every $f, g \in \cP^\mu(\MC{n})$. As we have done previously, Propositions~\ref{prop:isotypicabstract} and \ref{prop:a-Ta-Hinvariant} imply that the Toeplitz operators $T^{(\lambda)}_a$ and $(T^{(\lambda)}_a)^* = T^{(\lambda)}_{\overline{a}}$ preserve the subspace $\cP^\mu(\MC{n})$. 
	
	From the previous remarks we conclude that
	\begin{align*}
		T^{(\lambda)}_a|_{\cP^\mu(\MC{n})} &= T|_{\cP^\mu(\MC{n})} = T_A \\
		(T^{(\lambda)}_a)^*|_{\cP^\mu(\MC{n})} &= T^*|_{\cP^\mu(\MC{n})} = T_A^*,		
	\end{align*}
	and so that $T^{(\lambda)}_a$ is not normal.	
\end{proof}

\begin{remark}\label{rmk:LUnRUn-vs-End}
	The arguments used in the proof of Theorem~\ref{thm:notnormal-LUnRUn} are easily seen to imply some stronger properties. More precisely, if we consider the subspace of $\cA^2_\lambda(\DI{n})$ given by
	\[
		W = \bigoplus_{j=1}^k \cP^{\mu_j}(\MC{n})
	\]
	for a finite family of elements $\mu_1, \dots, \mu_k \in P(n)$, then for every operator $T \in \End_{\U(n)_L}(W) \subset \End_{\U(n)_L}(\cA^2_\lambda(\DI{n}))$ there exist symbols $a \in L^\infty(\DI{n})^{\U(n)_L}$ and $b \in L^\infty(\DI{n})^{\U(n)_R}$ such that
	\[
		T^{(\lambda)}_a|_W = T = T^{(\lambda)}_b|_W.
	\]
	In particular, for every $\mu \in P(n)$, the homomorphism of $C^*$-algebras given by
	\begin{align*}
		\cT^{(\lambda)}(L^\infty(\DI{n})^H) &\rightarrow \End_H(\cP^\mu(\MC{n})) \\
		T &\mapsto T|_{\cP^\mu(\MC{n})},
	\end{align*}
	is surjective when $H$ is taken to be either $\U(n)_L$ or $\U(n)_R$.
\end{remark}

\subsection{Commutative Banach algebras generated by Toeplitz operators}
\label{subsec:commutativeBanach}
In this subsection we will use both left and right $\U(n)$-invariant symbols to obtain commutative Banach algebras generated by Toeplitz operators that are not $C^*$-algebras. 

Our first examples are stated in the next result, which is a consequence of the existence of some special non-normal Toeplitz operators proved in Theorem~\ref{thm:notnormal-LUnRUn} (see also Proposition~\ref{prop:leftright-Un-invariant} and Remark~\ref{rmk:leftright-Un-invariant}).

\begin{theorem}\label{thm:noncomm-LUnRUn-oneoperator}
	Let us assume that $n \geq 2$. Then, for every $\lambda > 2n-1$ the following holds.
	\begin{enumerate}
		\item There exists $a \in L^\infty(\DI{n})$ that satisfies $a(Z) = a\big((Z^*Z)^\frac{1}{2}\big)$, for almost every $Z \in \DI{n}$, such that the unital Banach algebra generated by $T^{(\lambda)}_a$ is commutative and non-$C^*$. In particular, the unital $C^*$-algebra generated by $T^{(\lambda)}_a$ is not commutative.
		\item There exists $a \in L^\infty(\DI{n})$ that satisfies $a(Z) = a\big((ZZ^*)^\frac{1}{2}\big)$, for almost every $Z \in \DI{n}$, such that the unital Banach algebra generated by $T^{(\lambda)}_a$ is commutative and non-$C^*$. In particular, the unital $C^*$-algebra generated by $T^{(\lambda)}_a$ is not commutative.
	\end{enumerate}
\end{theorem}

We will also obtain commutative Banach algebras that contain the $C^*$-algebra $\cT^{(\lambda)}(L^\infty(\DI{n})^{\UU{n}})$ as well as Toeplitz operators with left and right $\U(n)$-invariant symbols. For this, we will use the next elementary auxiliary result.

\begin{lemma}\label{lem:vonNeumann-irreducible}
	Let $H$ be a Lie group with an irreducible unitary representation $\pi$ on a Hilbert space $\HH$. Then, the following properties are satisfied.
	\begin{enumerate}
		\item The von Neumann algebra generated by $\pi(H)$ is the algebra $\cB(\HH)$ of all bounded operators on $\HH$.
		\item If $\KK$ is a Hilbert space, and $\pi_L$ is the unitary representation of $H$ on $\HH \otimes \KK$ given by
		\[
			\pi_L(h)(u \otimes v) = (\pi(h)(u)) \otimes v
		\]
		for every $h \in H$, $u \in \HH$ and $v \in \KK$, then the von Neumann algebra generated by $\pi_L(H)$ consists of all the operators of the form $T \otimes I_\KK$ where $T \in \cB(\HH)$.
	\end{enumerate}
\end{lemma}
\begin{proof}
	By Schur's Lemma, the commutant $\pi(H)'$ is the algebra $\C I_{\HH}$. Hence, the von Neumann algebra generated by $\pi(H)$ is given by $\pi(H)'' = \cB(\HH)$. This proves~(1). Hence, (2) follows from (1) and the elementary properties of tensor products (see Corollary~1.5 in Chapter~IV of \cite{TakesakiOperatorI}).	
\end{proof}

We now obtain a centralizing result for the $C^*$-algebras considered in Theorem~\ref{thm:noncomm-LUnRUn}. It will be used to prove the existence of commutative Banach algebras as noted before. Nevertheless, this result is interesting by itself.

\begin{theorem}\label{thm:centralizing-LUnRUn}
	For every $\lambda \geq 2n -1$, the $C^*$-algebras $\cT^{(\lambda)}(L^\infty(\DI{n})^{\U(n)_L})$ and $\cT^{(\lambda)}(L^\infty(\DI{n})^{\U(n)_R})$ centralize each other. More precisely, we have
	\[
		ST = TS
	\]
	for every $S \in \cT^{(\lambda)}(L^\infty(\DI{n})^{\U(n)_L})$ and $T \in \cT^{(\lambda)}(L^\infty(\DI{n})^{\U(n)_R})$.
\end{theorem}
\begin{proof}
	By Proposition~\ref{prop:a-Ta-Hinvariant}(1) it is enough to prove the result for the von Neumann algebras $\End_{\U(n)_L}(\cA^2_\lambda(\DI{n}))$ and $\End_{\U(n)_R}(\cA^2_\lambda(\DI{n}))$. Since both of these algebras preserve the Hilbert direct sum from Theorem~\ref{thm:LUnRUn-isotypic}, it suffices to prove the result for the algebras $\End_{\U(n)_L}(\cP^\mu(\MC{n}))$ and $\End_{\U(n)_R}(\cP^\mu(\MC{n}))$, for every $\mu \in P(n)$. By Theorem~\ref{thm:UnUn-isotypic}, for every $\mu \in P(n)$ we have
	\[
		\cP^\mu(\MC{n}) \simeq W^{\mu *} \otimes W^\mu,
	\]
	as $\UU{n}$-modules. So the result finally reduces to showing that for every $\mu \in P(n)$, the von Neumann algebras $\End_{\U(n)_L}(W^{\mu *} \otimes W^\mu)$ and $\End_{\U(n)_R}(W^{\mu *} \otimes W^\mu)$ centralize each other.
	
	Since, both $W^\mu$ and $W^{\mu *}$ are irreducible $\U(n)$-modules,  Lemma~\ref{lem:vonNeumann-irreducible} shows that $\End_{\U(n)_L}(W^{\mu *} \otimes W^\mu)$ and $\End_{\U(n)_R}(W^{\mu *} \otimes W^\mu)$ consist precisely of maps of the form $T\otimes I_{W^\mu}$ and $I_{W^{\mu *}} \otimes S$, respectively, where $T \in \End(W^{\mu *})$ and $S \in \End(W^\mu)$. This clearly completes the proof.
\end{proof}

From the proof of Theorem~\ref{thm:centralizing-LUnRUn} we obtain the following immediate consequence.

\begin{corollary}\label{cor:vonNeumanncentralizing-LUnRUn}
	For every $\lambda > 2n-1$, the von Neumann algebras of intertwining operators $\End_{\U(n)_L}(\cA^2_\lambda(\DI{n}))$ and $\End_{\U(n)_R}(\cA^2_\lambda(\DI{n}))$ are the commutant of each other.
\end{corollary}

As a consequence of Theorems~\ref{thm:noncomm-LUnRUn-oneoperator} and \ref{thm:centralizing-LUnRUn} we obtain the following example of a commutative Banach algebra generated by Toeplitz operators.

\begin{corollary}\label{cor:noncomm-LUnRUn-twooperators}
	Let us assume that $n \geq 2$. For any given $\lambda > 2n-1$, let $a,b \in L^\infty(\DI{n})$ be symbols such that $a$ and $b$ satisfy (1) and (2) from Theorem~\ref{thm:noncomm-LUnRUn-oneoperator}, respectively. Then, the unital Banach algebra generated by $T^{(\lambda)}_a$ and $T^{(\lambda)}_b$ is commutative and non-$C^*$. In particular, the unital $C^*$-algebra generated by these two Toeplitz operators is not commutative.
\end{corollary}

Finally, we prove the existence of a commutative Banach algebra with a quite large set of generators which are Toeplitz operators. As before, for $n \geq 2$, the existence of the non-normal Toeplitz operators from the statement is guaranteed by Theorem~\ref{thm:notnormal-LUnRUn}. We recall that for a set of symbols $\cS$, we denote by $\cT^{(\lambda)}(\cS)$ the unital Banach algebra generated by the Toeplitz operators with symbols in $\cS$.

\begin{theorem}\label{thm:BanachCommNotC*}
	Assume that $n \geq 2$. For a given $\lambda > 2n-1$, choose symbols $a \in L^\infty(\DI{n})^{\U(n)_L}$ and $b \in L^\infty(\DI{n})^{\U(n)_R}$ such that the Toeplitz operators $T^{(\lambda)}_a$ and $T^{(\lambda)}_b$ are not normal, and let 
	\[
		\cS = L^\infty(\DI{n})^{\UU{n}} \cup \{a,b\}.
	\]
	Then, the Banach algebra $\cT^{(\lambda)}(\cS)$ is commutative and non-$C^*$. In particular, the unital $C^*$-algebra generated by $\cS$ is not commutative. 
\end{theorem}
\begin{proof}
	To prove the commutativity of the given Banach algebra $\cT^{(\lambda)}(\cS)$, it is enough to show that for every pair of symbols $\varphi, \psi \in \cS$ the Toeplitz operators $T^{(\lambda)}_\varphi$ and $T^{(\lambda)}_\psi$ commute with each other. If $\varphi, \psi$ both belong to $L^\infty(\DI{n})^{\UU{n}}$, then this follows from Theorem~\ref{thm:comm-UnUn}. If $\{\varphi, \psi\} = \{a,b\}$, then the claim follows from Theorem~\ref{thm:centralizing-LUnRUn}.
	
	Let us now assume that $\varphi \in L^\infty(\DI{n})^{\UU{n}}$ and that $\psi$ is either $a$ or $b$. In particular, $\psi$ is either left or right $\U(n)$-invariant. By Theorems~\ref{thm:comm-UnUn} and \ref{thm:noncomm-LUnRUn}, both Toeplitz operators $T^{(\lambda)}_\varphi$ and $T^{(\lambda)}_\psi$ preserve the Hilbert direct sum
	\[
		\cA^2_\lambda(\DI{n}) = \bigoplus_{\mu \in P(n)}
			\cP^\mu(\MC{n}),
	\]
	and the operator $T^{(\lambda)}_\varphi$ acts by a constant multiple of the identity on each term. Hence, the Toeplitz operators $T^{(\lambda)}_\varphi$ and $T^{(\lambda)}_\psi$ commute with each other. This proves that $\cT^{(\lambda)}(\cS)$ is commutative.
	
	Since $T^{(\lambda)}_a$ and $T^{(\lambda)}_b$ are not normal their adjoints do not belong to $\cT^{(\lambda)}(\cS)$. Hence, $\cT^{(\lambda)}(\cS)$ is not $C^*$ and so the unital $C^*$-algebra generated by $\cS$ is not commutative.
\end{proof}

\begin{remark}\label{rmk:otherBanachcommutative}
	For $n \geq 2$, and with our current notation, if we choose non-normal operators $S \in \cT^{(\lambda)}(L^\infty(\DI{n})^{\U(n)_L})$ and $T \in \cT^{(\lambda)}(L^\infty(\DI{n})^{\U(n)_R})$, then the unital Banach algebra generated by $\cT^{(\lambda)}(L^\infty(\DI{n})^{\UU{n}})$ together with $S,T$ is commutative, it is not $C^*$, and the unital $C^*$-algebra generated by such operators is not commutative. The proof is the same one used to obtain Theorem~\ref{thm:BanachCommNotC*}. This provides some more general commutative Banach non-$C^*$ algebras generated by Toeplitz operators.
	
	We also note that, if we take 
	\[
		\cS = L^\infty(\DI{n})^{\UU{n}} \cup \{\psi\},
	\]
	where $\psi \in \{a,b\}$ and $a,b$ are as in Theorem~\ref{thm:BanachCommNotC*}, then the Banach algebra $\cT^{(\lambda)}(\cS)$ satisfies the same conclusions from such theorem. This yields additional examples of commutative Banach non-$C^*$ algebras generated by Toeplitz operators. Of course, we can as well consider an even more general example given by the Banach algebra generated by $\cT^{(\lambda)}(L^\infty(\DI{n})^{\UU{n}})$ together with some non-normal operator in either $\cT^{(\lambda)}(L^\infty(\DI{n})^{\U(n)_L})$ or $\cT^{(\lambda)}(L^\infty(\DI{n})^{\U(n)_R})$.
\end{remark}


\begin{thebibliography}{XX}

\bibitem{AppuhamyLe2016} Appuhamy, Amila and Le, Trieu: \emph{Commutants of Toeplitz operators with separately radial polynomial symbols}. Complex Anal. Oper. Theory 10 (2016), no. 1, 1--12.

\bibitem{AxlerCuckovicRao2000} Axler, Sheldon; \v{C}u\v{c}kovi\'c, \v{Z}eljko and Rao, N. V.: \emph{Commutants of analytic Toeplitz operators on the Bergman space}. Proc. Amer. Math. Soc. 128 (2000), no. 7, 1951--1953.

\bibitem{BauerChoeKoo2015} Bauer, Wolfram; Choe, Boo Rim and Koo, Hyungwoon: \emph{Commuting Toeplitz operators with pluriharmonic symbols on the Fock space}. J. Funct. Anal. 268 (2015), no. 10, 3017--3060.

\bibitem{BHVRadial} Bauer, Wolfram; Herrera Yañez, Crispin and Vasilevski, Nikolai: \emph{Eigenvalue characterization of radial operators on weighted Bergman spaces over the unit ball}. Integral Equations Operator Theory 78 (2014), no. 2, 271--300.

\bibitem{BVquasir-quasih} Bauer, Wolfram and Vasilevski, Nikolai: \emph{On the structure of a commutative Banach algebra generated by Toeplitz operators with quasi-radial quasi-homogeneous symbols}. Integral Equations Operator Theory 74 (2012), no. 2, 199--231.

\bibitem{BVQuasiEllipticI}  Bauer, Wolfram and Vasilevski, Nikolai: \emph{On the structure of commutative Banach algebras generated by Toeplitz operators on the unit ball. Quasi-elliptic case. I: Generating subalgebras}. J. Funct. Anal. 265 (2013), no. 11, 2956--2990.

\bibitem{BtD} Br\"ocker, Theodor and tom Dieck, Tammo: Representations of compact Lie groups. Translated from the German manuscript. Corrected reprint of the 1985 translation. Graduate Texts in Mathematics, 98. Springer-Verlag, New York, 1995.

\bibitem{ChoeKooLee2004}  Choe, Boo Rim; Koo, Hyungwoon and Lee, Young Joo: \emph{Commuting Toeplitz operators on the polydisk}. Trans. Amer. Math. Soc. 356 (2004), no. 5, 1727--1749.

\bibitem{CuckovicLouhichi2008} \v{C}u\v{c}kovi\'c, \v{Z}eljko and Louhichi, Issam: \emph{Finite rank commutators and semicommutators of quasihomogeneous Toeplitz operators}. Complex Anal. Oper. Theory 2 (2008), no. 3, 429--439.

\bibitem{DOQJFA} Dawson, Matthew; \'Olafsson, Gestur and Quiroga-Barranco, Raul: \emph{Commuting   Toeplitz   operators   on   bounded symmetric   domains   and   multiplicity-free restrictions   of   holomorphic   discrete   series}. J. Funct. Anal. 268 (2015), no. 7, 1711--1732.

\bibitem{DQRadial} Dawson, Matthew and Quiroga-Barranco, Raul: \emph{Radial Toeplitz operators on the weighted Bergman spaces of Cartan domains}. Representation theory and harmonic analysis on symmetric spaces, 97–114, Contemp. Math., 714, Amer. Math. Soc., [Providence], RI, 2018.

\bibitem{EsmeralMaximenko2016} Esmeral, Kevin and Maximenko, Egor A.: \emph{Radial Toeplitz operators on the Fock space and square-root-slowly oscillating sequences}. Complex Anal. Oper. Theory 10 (2016), no. 7, 1655--1677.

\bibitem{GW} Goodman, Roe and Wallach, Nolan R.: Symmetry, representations, and invariants. Graduate Texts in Mathematics, 255. Springer, Dordrecht, 2009.

\bibitem{GKVRadial} Grudsky, S.; Karapetyants, A. and Vasilevski, N.: \emph{Toeplitz operators on the unit ball in $\C^n$ with radial symbols}. J. Operator Theory 49 (2003), no. 2, 325--346.

\bibitem{GKVRadialDisk} Grudsky, S.; Karapetyants, A. and Vasilevski, N.: \emph{Dynamics of properties of Toeplitz operators with radial symbols}. Integral Equations Operator Theory 50 (2004), no. 2, 217--253. 

\bibitem{GMVRadial} Grudsky, Sergei M.; Maximenko, Egor A. and Vasilevski, Nikolai L.: \emph{Radial Toeplitz operators on the unit ball and slowly oscillating sequences}. Commun. Math. Anal. 14 (2013), no. 2, 77--94.

\bibitem{GQVJFA} Grudsky, S.; Quiroga-Barranco, R. and Vasilevski, N.: \emph{Commutative $C^*$-algebras of Toeplitz operators and quantization on the unit disk}. J. Funct. Anal. 234 (2006), no. 1, 1–44. 

\bibitem{EnglisDensity} M. Engli\v{s}, Density of algebras generated by Toeplitz operators on Bergman spaces, Arkiv f\"or Matematik 30 (1992), No. 2, 227--243.

\bibitem{HuaHarmonic} Hua, L. K.: Harmonic analysis of functions of several complex variables in the classical domains. Translated from the Russian, which was a translation of the Chinese original, by Leo Ebner and Adam Kor\'anyi. With a foreword by M. I. Graev. Reprint of the 1963 edition. Translations of Mathematical Monographs, 6. American Mathematical Society, Providence, R.I., 1979.

\bibitem{KnappLieBeyond} Knapp, Anthony W. Lie groups beyond an introduction. Second edition. Progress in Mathematics, 140. Birkh\"auser Boston, Inc., Boston, MA, 2002.

\bibitem{KNVol1} Kobayashi, Shoshichi; Nomizu, Katsumi Foundations of differential geometry. Vol. I. Reprint of the 1963 original. Wiley Classics Library. A Wiley-Interscience Publication. John Wiley \& Sons, Inc., New York, 1996.

\bibitem{KorenblumZhu1995} Korenblum, Boris and Zhu, Ke He: \emph{An application of Tauberian theorems to Toeplitz operators}. J. Operator Theory 33 (1995), no. 2, 353--361.

\bibitem{Le2017} Le, Trieu: \emph{Commutants of separately radial Toeplitz operators in several variables}. J. Math. Anal. Appl. 453 (2017), no. 1, 48--63.

\bibitem{MokHermitian} Mok, Ngaiming: Metric rigidity theorems on Hermitian locally symmetric manifolds. Series in Pure Mathematics, 6. World Scientific Publishing Co., Inc., Teaneck, NJ, 1989.

\bibitem{QRadial} Quiroga-Barranco, Raul: \emph{Separately radial and radial Toeplitz operators on the unit ball and representation theory}. Bol. Soc. Mat. Mex. (3) 22 (2016), no. 2, 605--623.

\bibitem{QSpseudoconvex} Quiroga-Barranco, Raul and Sanchez-Nungaray, Armando: \emph{Toeplitz operators with quasi-homogeneous quasi-radial symbols on some weakly pseudoconvex domains}. Complex Anal. Oper. Theory 9 (2015), no. 5, 1111–-1134.

\bibitem{QVUnitBall1} Quiroga-Barranco, Raul and Vasilevski, Nikolai: \emph{Commutative $C^*$-algebras of Toeplitz operators on the unit ball. I. Bargmann-type transforms and spectral representations of Toeplitz operators.} Integral Equations Operator Theory 59 (2007), no. 3, 379--419.

\bibitem{QVUnitBall2} Quiroga-Barranco, Raul and Vasilevski, Nikolai: \emph{Commutative $C^*$-algebras of Toeplitz operators on the unit ball. II. Geometry of the level sets of symbols.} Integral Equations Operator Theory 60 (2008), no. 1, 89--132.

\bibitem{Rodriguez2020} Rodriguez, Miguel Angel Rodriguez: \emph{Banach algebras generated by Toeplitz operators with parabolic quasi-radial quasi-homogeneous symbols}. Bol. Soc. Mat. Mex. (3) 26 (2020), no. 3, 1243--1271.

\bibitem{SuarezRadial2008} Su\'arez, Daniel: \emph{The eigenvalues of limits of radial Toeplitz operators}. Bull. Lond. Math. Soc. 40 (2008), no. 4, 631--641.

\bibitem{Schmidt} Schmidt, Robert: \emph{\"Uber divergente Folgen and lineare Mittelbildungen}. Math. Z. 22 (1925), 89--152.

\bibitem{TakesakiOperatorI} Takesaki, Masamichi: Theory of operator algebras. I. Springer-Verlag, New York-Heidelberg, 1979.

\bibitem{UpmeierToepBook} Upmeier, Harald: Toeplitz operators and index theory in several complex variables. Operator Theory: Advances and Applications, 81. Birkh\"auser Verlag, Basel, 1996.

\bibitem{VasQuasiRadial} Vasilevski, Nikolai: \emph{Quasi-radial quasi-homogeneous symbols and commutative Banach algebras of Toeplitz operators}. Integral Equations Operator Theory 66 (2010), no. 1, 141--152.

\bibitem{VasPseudoH} Vasilevski, Nikolai: \emph{On Toeplitz operators with quasi-radial and pseudo-homogeneous symbols}, in: Harmonic Analysis, Partial Diﬀerential Equations, Banach Spaces, and Operator Theory, vol. 2, pp. 401--417, Springer Verlag, Heidelberg, 2017.

\bibitem{ZorbRadial2002} Zorboska, Nina: The \emph{Berezin transform and radial operators}. Proc. Amer. Math. Soc. 131 (2003), no. 3, 793--800.

\end{thebibliography}
\end{document}